%% file: main.tex
\documentclass{llncs}

\usepackage{algorithm}
\usepackage{algorithmic}
\usepackage{caption}
\usepackage{amsfonts}
\usepackage{amsmath}
\usepackage{amssymb}
\usepackage{hyperref}
\usepackage{listings}
\usepackage{makecell}
\usepackage[numbers]{natbib}
\usepackage{pgfplotstable}
\usepackage{relsize}
\hypersetup{colorlinks=true,linkcolor=black,citecolor=black,filecolor=black,urlcolor=black}
\pgfplotsset{compat=1.14}
\definecolor{palette1}{RGB}{228,26,28}
\definecolor{palette2}{RGB}{55,126,184}
\definecolor{palette3}{RGB}{77,175,74}
\definecolor{palette4}{RGB}{152,78,163}
\definecolor{palette5}{RGB}{255,127,0}

\begin{document}
	\input{frontmatter/frontmatter.tex}
	\input{introduction/introduction.tex}
	\input{literature_review/literature_review.tex}
	\input{background/background.tex}
	\input{linear_solvers/linear_solvers.tex}
	\input{computational_experiments/computational_experiments.tex}
	\input{concluding_remarks/concluding_remarks.tex}
	\input{acknowledgments/acknowledgments.tex}
	\input{bibliography/bibliography.tex}
	\input{appendix/appendix.tex}
\end{document}

%% file: frontmatter/frontmatter.tex
\input{frontmatter/authors.tex}
\input{frontmatter/title.tex}
\maketitle
\input{frontmatter/abstract.tex}
\input{frontmatter/keywords.tex}

%% file: frontmatter/authors.tex
\author{Byron Tasseff\inst{1,2} \and
        Carleton Coffrin\inst{1} \and
        Andreas W{\"a}chter\inst{3} \and
        Carl Laird\inst{4}}

\institute{Los Alamos National Laboratory, Los Alamos, New Mexico, USA \and
           University of Michigan, Ann Arbor, Michigan, USA \and
           Northwestern University, Evanston, Illinois, USA \and
           Sandia National Laboratories, Albuquerque, New Mexico, USA}

%% file: frontmatter/title.tex
\title{Exploring Benefits of Linear Solver Parallelism on Modern Nonlinear Optimization Applications}

%% file: frontmatter/abstract.tex
\begin{abstract}
	The advent of efficient interior point optimization methods has enabled the tractable solution of large-scale linear and nonlinear programming (NLP) problems.
	A prominent example of such a method is seen in \textsc{Ipopt}, a widely-used, open-source nonlinear optimization solver.
	Algorithmically, \textsc{Ipopt} depends on the use of a sparse symmetric indefinite linear system solver, which is heavily employed within the optimization of barrier subproblems.
	As such, the performance and reliability of \textsc{Ipopt} is dependent on the properties of the selected linear solver.
	Inspired by a trend in mathematical programming toward solving larger and more challenging NLPs, this work explores two core questions: first, how does the scalability of available linear solvers, many of which exhibit shared-memory parallelism, impact \textsc{Ipopt} performance; and second, does the best linear solver vary across NLP problem classes, including nonlinear network problems and problems constrained by partial differential equations?
	To better understand these properties, this paper first describes available open- and closed-source, serial and parallel linear solvers and the fundamental differences among them.
	Second, it introduces the coupling of a new open-source linear solver capable of heterogeneous parallelism over multi-core central processing units and graphics processing units.
	Third, it compares linear solvers using a variety of mathematical programming problems, including standard test problems for linear and nonlinear optimization, optimal power flow benchmarks, and scalable two- and three-dimensional partial differential equation and optimal control problems.
	Finally, linear solver recommendations are provided to maximize \textsc{Ipopt} performance across different application domains.
\end{abstract}

%% file: frontmatter/keywords.tex
\keywords{benchmarking \and graphics processing units \and interior point \and numerical linear algebra \and nonlinear programming \and parallel \and performance profiling}

%% file: introduction/introduction.tex
\section{Introduction}
\label{sec:introduction}
Since their early formation, beginning with \citet{frisch1955logarithmic} and \citet{fiacco1968nonlinear}, interior point (or barrier) methods have transformed the field of optimization.
Their low-degree polynomial complexity has made them especially attractive candidates for large-scale linear and convex optimization, as the number of iterations required for convergence is largely unaffected by problem dimension.
Interior point methods have proven to be competitive when applied to problems of small size (less than one million variables and constraints) and are without competition for problems of larger size \citep{GONDZIO2012587}.
As the demand to solve larger problems increases, presumably so will the utility of these methods.

The de facto open-source standard for interior point optimization is \textsc{Ipopt}, initially developed by \citet{wachter2006implementation} and implemented by Laird and W{\"a}chter.
\textsc{Ipopt} transforms inequality constraints into logartihmic barrier terms in the objective and solves a sequence of equality constrained problems for decreasing values of the barrier parameter.
For each barrier subproblem, a modified Newton algorithm with line search is used to obtain a solution to the problem's Karush-Kuhn-Tucker conditions.
Each step of the Newton algorithm requires potentially multiple solutions to perturbed sparse symmetric indefinite systems of linear equations.
Because of this inner algorithmic dependency, the performance of \textsc{Ipopt} depends critically on the properties and performance of the selected linear solver.
Note that \textsc{Ipopt}'s runtime is also dependent on the external computation of first and second derivatives, usually accomplished by the modeling language interface (e.g., \textsc{AMPL}, \textsc{GAMS}, \textsc{JuMP}, \textsc{Pyomo}).
The impact of improving the speed of these computations is not considered in this work.

The performance of a linear solver depends on the pivot sequence (or elimination order) selected prior to factorization.
Additionally, for symmetric indefinite solvers, this sequence may be modified during factorization to ensure numerical stability.
Due to the tradeoff between performance and robustness, a variety of techniques have been developed to perform this factorization.
Their differences in methodology, intended computational architecture, and software coupling to \textsc{Ipopt} can result in vastly different behaviors in the solution of nonlinear programs (NLPs).
This paper intends to provide a better understanding of these choices and their effects on the performance and robustness of \textsc{Ipopt}.

This paper provides three primary contributions, each aimed at effectively solving large and challenging NLPs via \textsc{Ipopt}.
The first is an understanding of the scalability of available linear solvers, many of which exhibit shared-memory parallelism.
The second is the coupling of a new open-source linear solver, \textsc{SPRAL}, capable of heterogeneous parallelism over multi-core central processing units (CPUs) and graphics processing units (GPUs).
The third is a comparison of all available linear solvers over a variety of common problem classes, including standard test problems for linear and nonlinear optimization, optimal power flow benchmarks, and two- and three-dimensional partial differential equation (PDE) and optimal control problems.
To summarize these contributions, the benefits of linear solver parallelism are assessed, and solver recommendations are provided to maximize \textsc{Ipopt}'s performance across different application domains.

The rest of this paper is organized as follows: Section \ref{sec:literature_review} discusses previous literature related to the benchmarking of interior point methods and the performance of \textsc{Ipopt}; Section \ref{sec:background} provides background information on interior point methods and the algorithm of \textsc{Ipopt}; Section \ref{sec:linear_solvers} compares the various solution methods used by the sparse symmetric indefinite linear solvers that are available within \textsc{Ipopt}; Section \ref{sec:computational_experiments} benchmarks these linear solvers using a variety of linear and nonlinear programming problems; and Section \ref{sec:concluding_remarks} concludes the paper.

%% file: literature_review/literature_review.tex
\section{Literature Review}
\label{sec:literature_review}
Historically, the benchmarking of optimization software has been an important component of solver development, frequently exposing numerical flaws and generally leading to software improvements \citep{dolan2002benchmarking}.
Throughout the various subfields of mathematical programming, a number of standard benchmarking libraries have emerged for solver comparison, including MIPLIB for mixed-integer linear programming and MINLPLIB for mixed-integer nonlinear programming \citep{koch2011miplib, vigerske2014minlplib}.
For general nonlinear programming, versions of the Constrained and Unconstrained Testing Environment (\textsc{CUTE}) of \citet{bongartz1995cute} have typically been used.
However, most benchmarking studies have either focused on accelerating the solution of specific problem types (e.g., \citet{rojas2015benchmarking, lyu2014benchmarking, blaszczyk2007object}) or broad comparisons among disparate solvers (e.g., \citet{dolan2004benchmarking}).
Few studies have aimed at benchmarking a single solver (e.g., \textsc{Ipopt}) with respect to the various options that can be employed by that solver.
As such, the goals of this study are twofold: first, to design and benchmark a modern set of numerical experiments relevant to large-scale nonlinear optimization; and second, to compare and provide useful guidance for selecting \textsc{Ipopt}'s most consequential solver option, the symmetric indefinite linear solver.

The study most similar to this one is that of \citet{schenk2007matching}.
Here, the integration of the \textsc{PARDISO} linear solver within \textsc{Ipopt} is introduced.
Several preprocessing and pivoting strategies are compared using \textsc{CUTE}, the Constrained Optimization Problem Set (\textsc{COPS}) of \citet{dolan2004benchmarking}, and a number of large PDE-constrained optimal control problems.
To demonstrate its effectiveness, use of \textsc{PARDISO} is compared against two other linear solvers.
For the problems considered, \textsc{PARDISO} performs favorably.
However, compared to this paper, their study focuses on only three linear solvers and limited problem sets.

Another similar study is that of \citet{hogg2013effects}.
Here, they explore various methods for scaling symmetric indefinite matrices within \textsc{Ipopt} while using the \textsc{HSL\_MA97} linear solver.
Using a subset of large experiments from the Constrained and Unconstrained Testing Environment (Revisited; \textsc{CUTEr}) of \citet{gould2003cuter}, they compare these methods with respect to their impacts on the performance and reliability of \textsc{Ipopt}.
In general, they find the benefits of scaling are problem dependent and thus propose a number of heuristics for switching among methods during the solution process.
Compared to this paper, their study considers only one linear solver, whereas this study considers a variety of available linear solvers.
Additionally, only a small set of problems is examined in their study, whereas this study considers larger and more comprehensive sets.

In general, \textsc{Ipopt} assumes the use of direct (noniterative) factorization methods in the solution of symmetric indefinite systems.
A number of studies have proposed alternative strategies that enable the use of iterative methods.
\citet{schenk2008inertia} propose an iterative technique for \textsc{PARDISO} that performs well on large-scale NLPs.
\citet{curtis2010interior} propose an interior point algorithm that leverages the efficiency of inexact step computations with iterative linear solvers.
\citet{curtis2012note} provide an implementation of this algorithm within \textsc{Ipopt} that uses \textsc{PARDISO}, finding good performance on PDE-constrained problems.
\citet{grote2014inexact} extend this method and find success on nonconvex PDE-constrained problems.
However, as of this writing, iterative methods available in \textsc{Ipopt} are largely experimental.
Thus, this study only considers linear solvers that use direct methods, based on the algorithmic requirements described in Section \ref{sec:background}.

%% file: background/background.tex
\section{Background}
\label{sec:background}
Consider a general mathematical programming problem of the form
\begin{subequations}\begin{alignat}{2}
	& \underset{x \in \mathbb{R}^{n}}{\text{minimize}}
	& & ~~~ f(x) \\
	& \text{subject to}
	& & ~~~ c(x) = 0 \\
	& & & ~~~ x \geq 0,
\end{alignat}\label{eqn:gen-opt}\end{subequations}
where $f : \mathbb{R}^{n} \to \mathbb{R}$ is the objective function and $c : \mathbb{R}^{n} \to \mathbb{R}^{m}$ are constraints, all of which are assumed to be twice continuously differentiable.
Using duality theory, the first-order optimality conditions for a problem of this form are
\begin{subequations}\begin{align}
	\nabla f(x) + \nabla c(x) y - z &= 0 \\
	c(x) &= 0 \\
	X Z e &= 0 \label{eqn:gen-complementarity} \\
	x, z &\geq 0,
\end{align}\label{eqn:gen-fonc}\end{subequations}
where $y \in \mathbb{R}^{m}$ and $z \in \mathbb{R}^{n}$ are the Lagrangian dual multipliers for the equality and bound constraints, respectively.
Additionally, $X = \textrm{Diag}(x)$, $Z = \textrm{Diag}(z)$, and $e = (1, 1, \dots, 1)^{\intercal}$.
Equation \eqref{eqn:gen-complementarity} is the complementarity condition and implies that one of the two variables, $x_{j}$ and $z_{j}$, must be equal to zero at optimality.
The method by which the complementarity condition is handled determines the classification of the corresponding solution algorithm, which will typically fall into one of two categories: \emph{active set methods} or \emph{interior point methods}.

Active set methods operate by selecting subsets of indices $\mathcal{B} \subset \{1, 2, \dots, n\}$ such that $x_{j} \neq 0, ~ \forall j \in \mathcal{B}$.
An example of such a method is the primal simplex algorithm for linear programming, which sequentially separates variables into active (basic) and inactive (nonbasic) sets.
Interior point methods, on the other hand, replace $x_{j} z_{j} = 0$ with $x_{j} z_{j} = \mu$, where $\mu$ is a parameter that is gradually driven toward zero.
Thus, with interior point methods, the partitioning of variables into active and inactive sets is revealed as optimization progresses.
To accomplish this, interior point methods usually exploit three tools: first, logarithmic barrier functions to transform inequality constraints into objective function terms; second, duality theory to derive first-order optimality conditions; and third, Newton-type algorithms to solve systems of first-order conditions \citep{GONDZIO2012587}.

\textsc{Ipopt} is an interior point algorithm that considers a variant of Problem \eqref{eqn:gen-opt},
\begin{subequations}\begin{alignat}{2}
	& \underset{x \in \mathbb{R}^{n}}{\text{minimize}}
	& & ~~~ f(x) \\
	& \text{subject to}
	& & ~~~ g^{L} \leq g(x) \leq g^{U} \\
	& & & ~~~ x^{L} \leq x \leq x^{U},\end{alignat}\label{ipopt-opt}\end{subequations}
where $x^{L}, x^{U} \in \mathbb{R}^{n}$ are lower and upper bounds on optimization variables and $g : \mathbb{R}^{n} \to \mathbb{R}^{m}$ are constraints with lower and upper bounds, $g^{L}$ and $g^{U}$, respectively.
\textsc{Ipopt} first transforms problems of this form into a form similar to Problem $\eqref{eqn:gen-opt}$, where slack variables are introduced such that $g_{i}(x) - s_{i} = 0$ and $g_{i}^{L} \leq s_{i} \leq g_{i}^{U}$.
For notational ease, the remainder of this section assumes all variables (original and slack) only have lower bounds of zero.
The more general case is described by \citet{wachter2006implementation}.
To find a locally optimal solution to the problem, \textsc{Ipopt} then solves a converging sequence of barrier subproblems of the form
\begin{subequations}
\begin{alignat}{2}
	& \underset{x \in \mathbb{R}^{n}}{\text{minimize}}
	& & ~~~ \phi_{\mu}(x) := f(x) - \mu \sum_{i = 1}^{n} \ln(x_{i}) \\
	& \text{subject to}
	& & ~~~ c(x) = 0.
\end{alignat} \label{eqn:ipopt-barrier}
\end{subequations}

For each barrier subproblem, the first-order optimality conditions are thus
\begin{subequations}\begin{align}
	\nabla f(x) + \nabla c(x) y - z &= 0 \label{eqn:barrier-condition-1} \\
	c(x) &= 0 \label{eqn:barrier-condition-2} \\
	X Z e - \mu e &= 0 \label{eqn:barrier-condition-3} \\
	x, z &\geq 0.
\end{align}\label{eqn:barrier-conditions}\end{subequations}
Note that when $\mu = 0$, these are equivalent to the conditions described by System \eqref{eqn:gen-fonc}.
Under certain assumptions, it can be shown that the optimal solutions $x^{*}(\mu)$ of Problem \eqref{eqn:ipopt-barrier} converge to an optimal solution of Problem \eqref{eqn:gen-opt} as $\mu$ approaches zero \citep{byrd1997local, gould2001superlinear}.
Thus, within \textsc{Ipopt}, $\mu$ is monotonically decreased until a point is obtained that satisfies first-order optimality conditions (to some predefined tolerance).
However, it is important to note that maximizers and saddle points may also satisfy these conditions.
Because of this, \textsc{Ipopt} does not guarantee convergence to a local minimizer.
However, to address this, a regularization technique partially described below encourages \textsc{Ipopt} to avoid these points \citep{wachter2009short}.

For a fixed $\mu = \bar{\mu}$, the corresponding barrier problem is solved using a Newton-type algorithm.
This algorithm generates a converging sequence of iterates $(x_{k}, y_{k}, z_{k})$, with each iterate satisfying $x_{k}, z_{k} > 0$.
For each iterate, the step $(\Delta x_{k}, \Delta y_{k}, \Delta z_{k})$ can be obtained by solving for the roots of the first-order Taylor expansion of Equations \eqref{eqn:barrier-condition-1} through \eqref{eqn:barrier-condition-3} at $(x_{k}, y_{k}, z_{k})$, that is,
\begin{equation}
	\begin{pmatrix}
		W_{k} & \nabla c(x_{k}) & -I \\
		\nabla c(x_{k})^{\intercal} & 0 & 0 \\
		Z_{k} & 0 & X_{k}
	\end{pmatrix}
	\begin{pmatrix}
		\Delta x_{k} \\
		\Delta y_{k} \\
		\Delta z_{k}
	\end{pmatrix} =
	-\begin{pmatrix}
		\nabla f(x_{k}) + \nabla c(x_{k}) y_{k} - z_{k} \\
		c(x_{k}) \\
		X_{k} Z_{k} e - \bar{\mu} e
	\end{pmatrix},
\end{equation}
where $W_{k}$ is the Hessian of the Lagrangian function with respect to $x_{k}$, that is,
\begin{equation}
	W_{k} = \nabla^{2} f(x_{k}) + \sum_{j = 1}^{m} y_{k}^{j} \nabla^{2} c^{j}(x_{k}).
\end{equation}
\textsc{Ipopt} elects to solve this system by instead solving the smaller symmetric system
\begin{equation}
	\begin{pmatrix}
		W_{k} + X_{k}^{-1} Z_{k} & \nabla c(x_{k}) \\
		\nabla c(x_{k})^{\intercal} & 0
	\end{pmatrix}
	\begin{pmatrix}
		\Delta x_{k} \\
		\Delta y_{k}
	\end{pmatrix} =
	-\begin{pmatrix}
		\nabla \phi_{\bar{\mu}}(x_{k}) + \nabla c(x_{k}) y_{k} \\
		c(x_{k})
	\end{pmatrix} \label{eqn:iteration-matrix-orig}.
\end{equation}
From the solution to this system, $\Delta z_{k}$ can then be reconstructed via
\begin{equation}
	\Delta z_{k} := \bar{\mu} X_{k}^{-1} e - z_{k} - X_{k}^{-1} Z_{k} \Delta x_{k}.
\end{equation}
With $\alpha_{k}^{x, y, z} \in (0, 1]$ serving as step sizes, the iterate $(x_{k+1}, y_{k+1}, z_{k+1})$ is then
\begin{subequations}
\begin{align}
	x_{k+1} &= x_{k} + \alpha_{k}^{x} \Delta x_{k} \\
	y_{k+1} &= y_{k} + \alpha_{k}^{y} \Delta y_{k} \\
	z_{k+1} &= z_{k} + \alpha_{k}^{z} \Delta z_{k}.
\end{align}
\end{subequations}

\textsc{Ipopt} encourages convergence to a local minimizer through its methodology for determining the step size $(\alpha_{k}^{x}, \alpha_{k}^{y}, \alpha_{k}^{z})$, as well as its methodology for making second-order corrections that improve trial points in the Newton-type algorithm.
Specifically, a line search variant of the filter method described by \citet{fletcher2002nonlinear} is used to compute a step size $(\alpha_{k}^{x}, \alpha_{k}^{y}, \alpha_{k}^{z})$ that attempts to ensure either a decrease in the objective of the barrier subproblem or a decrease in the overall constraint violation.
The details of this line search method are omitted here for brevity, but the reader is referred to \citet{wachter2006implementation} for elaboration.
It is suffice to say that in order to make this guarantee, the iteration matrix of Equation \eqref{eqn:iteration-matrix-orig} must be nonsingular, and $(W_{k} + X_{k}^{-1} Z_{k}) v$, where $\nabla c(x_{k})^{\intercal} v = 0$, must be positive definite.
To satisfy these two properties, the iteration matrix must have exactly $n$ positive eigenvalues, $m$ negative eigenvalues, and no zero eigenvalues \citep{nocedal2006numerical}.
As such, the system may have to be perturbed as
\begin{equation}
	\begin{pmatrix}
		W_{k} + X_{k}^{-1} Z_{k} + \delta_{w} I & \nabla c(x_{k}) \\
		\nabla c(x_{k})^{\intercal} & -\delta_{c} I
	\end{pmatrix}
	\begin{pmatrix}
		\Delta x_{k} \\
		\Delta y_{k}
	\end{pmatrix} =
	-\begin{pmatrix}
		\Delta \phi_{\bar{\mu}}(x_{k}) + \nabla c(x_{k}) y_{k} \\
		c(x_{k})
	\end{pmatrix}, \label{eqn:ipopt-iterate-solve}
\end{equation}
where $\delta_{w}, \delta_{c} \geq 0$.
The values of $\delta_{w}$ and $\delta_{c}$ are thus varied heuristically until the eigenvalues are as desired.
(Note that there always exist $\delta_{w}, \delta_{c}$ to ensure this.)

A similar method is used to perform second-order corrections to iterates $(x_{k}, y_{k}, z_{k})$.
These corrections aim at reducing the infeasibility of the current primal iterate, $x_{k}$.
To accomplish this, $c(x_{k})$ on the right-hand side of Equation \eqref{eqn:ipopt-iterate-solve} is replaced with some perturbed vector $c_{k}^{\textrm{soc}}$.
Similar to the inertia correction, this perturbation process is repeated until the iterate satisfies acceptance criteria.

Because of these two inner algorithmic dependencies, \textsc{Ipopt} depends critically on a linear solver that carries with it two features: first, it should compute the numbers of positive, negative, and zero eigenvalues (``inertia'') of symmetric indefinite matrices; and second, it should efficiently reuse information from previous solutions to Equation \eqref{eqn:ipopt-iterate-solve} when solving for different right-hand sides.
The former property is a restrictive requirement, as inertia is typically (but not always) calculated as a byproduct of symmetric indefinite linear solvers that use direct factorization techniques.
This study thus focuses on the highly specific but computationally important investigation of the effect of using various symmetric indefinite direct linear solvers on the performance and robustness of \textsc{Ipopt}.

%% file: linear_solvers/linear_solvers.tex
\section{Solvers for Sparse Symmetric Indefinite Linear Systems}
\label{sec:linear_solvers}

Consider the system of equations presented by Equation \eqref{eqn:ipopt-iterate-solve}, rewritten as
\begin{equation}
	A x = b, \label{eqn:axeqb}
\end{equation}
where $A$ is sparse, symmetric, and indefinite.
Direct methods broadly divide the solution process into four phases: ordering, analysis, factorization, and solution.
The order phase computes a pivot sequence (ordering) intended to minimize the number of nonzeros (``fill-in'') in matrix factors.
The analyze phase uses this ordering and the sparsity pattern of $A$ to construct data structures for factorization.
In this phase, a scaling of the matrix can also be computed to improve numerical performance.
The factorize phase prepares the decomposition
\begin{equation}
	A = P A P^{\intercal} = L D L^{\intercal},
\end{equation}
where $P$ is a permutation matrix holding the pivot order, $L$ is a unit lower triangular matrix, and $D$ is a block diagonal matrix with blocks of order one or two.
The solve phase begins by first performing forward substitution, where
\begin{equation}
	L y = P b
\end{equation}
is solved for $y$.
Next, the block diagonal system
\begin{equation}
	D z = y
\end{equation}
is solved for $z$, followed by backward substitution, where
\begin{equation}
	L^{\intercal} P^{-\intercal} x = z
\end{equation}
is solved for $x$.
Note that some solvers are capable of using information from prior phases to efficiently solve Equation \eqref{eqn:axeqb} for multiple right-hand sides \citep{hogg2010high}.

The primary differences among symmetric indefinite solvers are in the factorize phase, which is also the most computationally expensive.
The linear solvers considered herein are thus separated into two categories based on their factorization technique: \emph{multifrontal methods} and \emph{supernodal methods}.
The \emph{frontal method} predates both and involves the use of a single dense ``frontal'' submatrix that holds the portion of the sparse matrix actively being factorized, with rows and columns entering and exiting during factorization.
The multifrontal method divides the original sparse matrix into \emph{several} frontal matrices, related to one another via a tree that represents dependencies between partial factorizations.
This allows each frontal matrix to exploit dense matrix operations and encourages parallelism.
The supernodal method combines adjacent columns with identical nonzero patterns from the sparse matrix into separate ``supernodes.''
Like multifrontal methods, supernodal methods also exploit dense matrix operations.
However, supernodal methods are typically more memory efficient \citep{davis2016survey}.

A detailed description of the theoretical and numerical evolution of sparse direct solvers is provided by \citet{davis2016survey}.
In this study, most of these details are omitted for brevity.
A summary of available sparse symmetric indefinite linear solvers within \textsc{Ipopt} can be found in Table \ref{tab:linear-solver-summary}.
Here, the factorization method, parallel capabilities, and license restrictions for each solver are described.
Solver properties relevant to this study are elaborated upon in Sections \ref{sec:multifrontal} and \ref{sec:supernodal}.

\begin{table}[t]
	\begin{center}
		\begin{tabular}{|c|c|c|c|c|}
			\hline \textbf{Solver} & \textbf{Method}  & \textbf{Parallel CPU} & \textbf{Parallel GPU} & \textbf{License Restrictions} \\ \hline
			\textsc{MA27}        & Multifrontal & No            & No           & \makecell{Free to all \\ Redistribution prohibited} \\ \hline
			\textsc{MA57}        & Multifrontal & Threaded BLAS & No           & Free to academics \\ \hline
			\textsc{HSL\_MA77}   & Multifrontal & Limited & No           & Free to academics \\ \hline
			\textsc{HSL\_MA86}   & Supernodal   & \makecell{Multi-core \\ Threaded BLAS} & No           & Free to academics \\ \hline
			\textsc{HSL\_MA97}   & Multifrontal & \makecell{Multi-core \\ Threaded BLAS}    & No           & Free to academics \\ \hline
			\textsc{MUMPS}       & Multifrontal & \makecell{Multi-core \\ Threaded BLAS}    & No           & Free to all \\ \hline
			\textsc{PARDISO}     & Supernodal   & Multi-core    & No           & Academic/corporate license \\ \hline
			\textsc{SPRAL}       & Multifrontal & \makecell{Multi-core \\ Threaded BLAS}    & Yes          & Free to all \\ \hline
			\textsc{WSMP}        & Multifrontal & Multi-core    & No           & Trial/purchased license \\ \hline
		\end{tabular}
		\caption{Summary of linear solvers available within \textsc{Ipopt}, including their factorization methodologies, parallel capabilities, and license restrictions.}
		\label{tab:linear-solver-summary}
	\end{center}
	\vspace{-3em}
\end{table}

\subsection{Multifrontal methods}
\label{sec:multifrontal}
\paragraph{MA27}
\label{subsec:ma27}
\citet{ma27b} introduce \textsc{MA27}, the first multifrontal solver ever developed.
Specifically, \textsc{MA27} focuses on the solution of symmetric systems, with a primary focus on the stable solution of symmetric indefinite systems.
This is achieved by allowing the pivot order to be modified if a pivot does not satisfy stability criteria.
The solver also allows for the solution of multiple right-hand sides \citep{gould2004numerical}.
Although phases of the solver are vectorized, they are not parallel.

\paragraph{MA57}
\label{subsec:ma57}
\citet{ma57} introduces \textsc{MA57}, designed to supersede \textsc{MA27} for the solution of symmetric indefinite systems.
One of the primary reasons for this supersession is that, since \textsc{MA27} was developed in the 1980s, it does not exploit Basic Linear Algebra Subprogram (BLAS) calls of Level 2 (matrix-vector operations) or Level 3 (matrix-matrix operations).
To increase portability, these are used throughout \textsc{MA57}.
The solver also implements the ability to restart and continue factorization if memory becomes limited \citep{gould2004numerical}.
Although \textsc{MA57} is sequential, it can exploit the parallelism of multithreaded BLAS calls.
\textsc{MA57} is developed for Fortran 77 and is thus different from \textsc{HSL\_MA57}, which is a Fortran 95 encapsulation.

\paragraph{HSL\_MA77}
\label{subsec:ma77}
\citet{ma77a, ma77b} describe the design and implementation of \textsc{HSL\_MA77}, an out-of-core (or external memory) solver for symmetric indefinite systems.
The solver implements a custom virtual memory scheme that allows data to be stored outside of main memory.
In turn, this allows for the solution of systems that could not normally reside in main memory.
Unless the solve phase is called repeatedly, the overhead of the method is not prohibitive.
The solver also includes an in-core method for when out-of-core methods are not required.
Similar to \textsc{MA57}, \textsc{HSL\_MA77} is serial but can exploit multithreaded BLAS.

\paragraph{HSL\_MA97}
\label{subsec:ma97}
\citet{ma97} develop \textsc{HSL\_MA97} as a means to explore modern multifrontal methods.
Their motivation is a number of deficiencies in \textsc{MA57} and \textsc{HSL\_MA77}, including their lack of modularity and parallelism.
Specifically, \textsc{HSL\_MA97} incorporates a shared-memory parallel method implemented using OpenMP.
Compared to other parallel sparse direct solvers, \textsc{HSL\_MA97} also focuses on the notion of ``bit-compatibility,'' or the bit-for-bit reproduction of results, regardless of the number of parallel threads used in the solution process.

\paragraph{MUMPS}
\label{subsec:mumps}
\citet{amestoy2000mumps} develop the Multifrontal Massively Parallel Solver (\textsc{MUMPS}).
This solver supports a number of features, including shared-memory parallelism via multithreaded BLAS and OpenMP.
Out-of-core capabilities similar to \textsc{HSL\_MA77} are also supported.
However, the primary design of \textsc{MUMPS} is geared toward distributed computation using MPI.
Note, though, that the focus of this study is on the comparison of linear solvers in a shared-memory environment, so the performance of \textsc{MUMPS} in \textsc{Ipopt} should be properly caveated.

\paragraph{WSMP}
\label{subsec:wsmp}
\citet{gupta2000wsmp} describes the direct symmetric solution method used by the Watson Sparse Matrix Package (\textsc{WSMP}).
Similar to \textsc{MUMPS}, \textsc{WSMP} is capable of shared-memory and distributed parallelism, although its indefinite solver is only amenable to the former.
\textsc{WSMP} performs its own parallelization and assumes the use of serial BLAS calls.
A fully serial version is also available.

\paragraph{SPRAL}
\citet{Hogg:2016:SSI:2888419.2756548} and \citet{duff2018new} describe the implementation of a sparse symmetric indefinite direct solver (SSIDS) within the the Sparse Parallel Robust Algorithms Library (SPRAL).
The motivation for the SPRAL SSIDS (hereafter abbreviated as \textsc{SPRAL}) is the limited parallelism of conventional pivoting techniques in multifrontal indefinite factorization.
To address this, they introduce ``a posteriori threshold pivoting,'' which more efficiently exploits parallel architectures while preserving stability.
The solver was initially implemented for parallelization over NVIDIA GPUs but now enables heterogeneous parallelization over multi-core CPUs, as well.
Similar to \textsc{HSL\_MA97}, when used with a bit-compatible BLAS library, \textsc{SPRAL} guarantees bit-compatibility of results.

Prior to this study, \textsc{SPRAL} was not available within \textsc{Ipopt}.
This is surprising due to its stability-preserving features and heterogeneous parallelization capabilities.
The library interface to the solver is also remarkably similar to that of \textsc{HSL\_MA97}, making integration with \textsc{Ipopt} straightforward.
For these reasons, aside from the benchmarking of currently available linear solvers, this paper also introduces the integration of \textsc{SPRAL} within \textsc{Ipopt} and studies its performance on multi-core CPUs and GPUs.
For the most part, the integration follows that of \textsc{HSL\_MA97} and retains similar default numerical parameterizations used across \textsc{HSL} solvers.
Additionally, parameters that control the various scaling, ordering, and parallel methods within \textsc{SPRAL} are provided as user options within \textsc{Ipopt}.

\subsection{Supernodal methods}
\label{sec:supernodal}
\paragraph{PARDISO}
\label{subsec:pardiso}
\citet{schenk2001pardiso} introduce \textsc{PARDISO}, a shared-memory parallel supernodal method.
For symmetric indefinite systems, it uses a technique based on a symmetric weighted matching algorithm to determine a pivot order that remains static during factorization.
Note that this contrasts with the threshold pivoting approach employed by most multifrontal solvers, where pivots are delayed if certain stability criteria are not satisfied.
The overall strategy of \textsc{PARDISO} thus ensures that reasonable pivots are predicted within supernodes at the cost of numerical stability, which generally gives the solver a preference toward performance over robustness.
The efficacy of this approach in \textsc{Ipopt} has been studied by \citet{schenk2007matching}.
Since its initial release, substantial improvements have continually been made, for example by \citet{schenk2002two, schenk2006fast}.
It has also been successfully employed and extended over a variety of applications, for example by \citet{pardiso-6.0a}, \citet{pardiso-6.0b}, and \citet{pardiso-6.0c}.

\paragraph{HSL\_MA86}
\label{subsec:ma86}
\citet{ma86} introduce \textsc{HSL\_MA86}, which uses a directed acyclic graph- (DAG-) based algorithm that facilitates fine-grained multi-core parallelism during factorization.
Here, as in other supernodal methods, the DAG describes dependencies between partial factorization tasks.
In the indefinite case, if a pivot candidate is delayed due to instability, matrix columns may be moved to different block columns and/or DAG nodes.
If many pivots are delayed, this can substantially reduce the level of parallelism employed by the solver.
Additionally, for hard indefinite matrices, \textsc{HSL\_MA86} includes only one ordering technique, \textsc{HSL\_MC68}, which may result in more delays during factorization than when used with a more robust ordering technique (e.g., matching-based ordering) \citep{duff2018new}.

%% file: computational_experiments/computational_experiments.tex
\section{Computational Experiments}
\label{sec:computational_experiments}
In this section, numerical results are presented for the solution of various linear and nonlinear programs using \textsc{Ipopt} and the nine available linear solvers.
In Section \ref{subsec:calibration}, parameters for linear solvers with parallel capabilities are calibrated to enable defensible comparisons in subsequent experiments.
In Section \ref{subsec:cutest}, performance is compared on problems from the Constrained and Unconstrained Testing Environment (safe threads; \textsc{CUTEst}) \citep{Gould:2015:CCU:2746412.2746437}.
In Section \ref{subsec:opf}, various alternating current optimal power flow (AC-OPF) formulations are considered with version 19.05 of the ``Power Grid Lib - Optimal Power Flow'' of \citet{pglibopf}.
In Section \ref{subsec:2d-pde}, scalable versions of problems by \citet{maurer2000optimization, maurer2001optimization} and \citet{mittelmann2001sufficient} are considered, all of which are constrained by a set of two-dimensional (2D) PDEs.
Finally, in Section \ref{subsec:3d-pde}, three-dimensional (3D) variants of the boundary control problems by \citet{maurer2000optimization} are considered.

Experiments were performed on Los Alamos National Laboratory's Darwin computing cluster.
All CPU-based experiments were executed on a node containing two Intel Xeon E5-2695 v4 processors, each with 18 cores @2.10 GHz, and 125 GB of memory.
All GPU-based experiments were executed on a node containing two Intel Xeon E5-2660 v3 processors, each with 10 cores @2.60 GHz, 125 GB of memory, and a single NVIDIA Tesla K40M GPU containing 2880 CUDA cores @745 MHz and 12 GB of memory.
For GPU-based experiments, \textsc{SPRAL}'s \texttt{gpu\_perf\_coeff} and \texttt{min\_gpu\_work} parameters were set to 999.9 and zero, respectively, in an attempt to induce maximal GPU usage for comparison.

For \textsc{SPRAL}, \textsc{METIS}-based ordering and \textsc{MC64} scaling were used for all \textsc{CUTEst} and OPF experiments, while matching-based ordering and scaling were used for PDE-constrained experiments.
Additionally, a modified \textsc{Ipopt} interface to \textsc{PARDISO}, provided by Andreas W{\"a}chter, was used.
Otherwise, all linear solvers used their default \textsc{Ipopt} settings.
The \texttt{acceptable\_tol} option in \textsc{Ipopt} was given a value of $10^{-8}$ (equivalent to the default \texttt{tol}).
Each experiment was given a wall clock time limit of four hours and a maximum \textsc{Ipopt} iteration limit of $9999$.
Experiments that resulted in optimal or ``acceptable'' solution statuses within these limits were considered ``solved.''
Accounting for early completion, the full set of experiments required over $82$ days of (wall clock) compute time.

Compilation required the Intel C\textsc{++} and Fortran 19.0.4 compiler suite, used only for the Intel Math Kernel Library (MKL) BLAS and LAPACK libraries contained therein.
All software was compiled using the GNU Compiler Collection (GCC; 7.2.0).
For GPU libraries and compilation, NVIDIA's CUDA Toolkit (10.1) was used.
Additionally, GNU Libtool (2.4.6), Portable Hardware Locality (hwloc; 1.11.12), and METIS (4.0.3) were self-compiled.
For HSL-based linear solvers (i.e., those with prefixes \textsc{HSL} and \textsc{MA}), the 2014-01-10 ``Coin-HSL Full (Stable)'' library was used.
For \textsc{MUMPS}, version 4.10.0 (whose installation script is provided with \textsc{Ipopt}) was used.
For \textsc{PARDISO}, the GNU 7.2.0 variant of the version 6.0.0 library was used.
For \textsc{WSMP}, the GNU variant of the 18.06.02 library was used. 
For \textsc{SPRAL}, the 2016-09-23 source release was self-compiled, with minor fixes and modifications made for integration with \textsc{Ipopt}.
\textsc{PARDISO} and \textsc{WSMP} assumed the use of sequential BLAS libraries.
As such, the \textsc{Ipopt} build corresponding to these solvers was compiled using sequential MKL flags.
Additional compilation information is provided in the appendix.

\subsection{Calibration of Parallel Linear Solvers}
\label{subsec:calibration}
As described in Table \ref{tab:linear-solver-summary}, the parallelism of linear solvers considered throughout this study is variable.
Although \textsc{SPRAL} has internal methods based on the Portable Hardware Locality library for ascertaining parallel resources and distributing among them, for other solvers, these capabilities are limited.
Additionally, because of the overhead of inter-socket communication on multi-CPU systems, using the default (maximum) number of threads for parallel processes could be a suboptimal choice for some linear solvers.
Thus, to promote fairer comparisons among linear solvers, this study elects to calibrate the number of parallel threads used by select solvers for subsequent experiments and analysis.
For \textsc{HSL\_MA86}, \textsc{HSL\_MA97}, \textsc{MUMPS}, and \textsc{PARDISO}, this entails calibrating the environment variable \texttt{OMP\_NUM\_THREADS}, and for \textsc{WSMP}, \texttt{WSMP\_NUM\_THREADS}.
For \textsc{SPRAL}, calibration is omitted because of its internal capabilities.
For remaining solvers, parallelism is limited, so calibration is also not performed.

To perform this calibration, a three-dimensional PDE-constrained boundary control problem similar to that presented in Example 5.1 of \citet{maurer2000optimization} is considered.
In this problem, the domain is $\Omega = (0, 1) \times (0, 1) \times (0, 1)$, and the goal is to compute the optimal boundary control $u(x)$ and state $y(x)$ with respect to $x = (x_{1}, x_{2}, x_{3})$.
The example is implemented in \textsc{Ipopt}'s \texttt{examples/ScalableProblems/solve\_problem} program as \texttt{MBndryCntrl\_3D}.
In the calibration herein, the scaling dimension $N = 48$ is used, which results in a nonlinear program containing $124,416$ variables ($13,824$ of which have lower and upper bounds and $110,592$ of which have only upper bounds), $110,592$ equality constraints, and $774,144$ nonzeros in the Jacobian matrix.
The problem is sufficiently large, results in substantial fill-in, and typically requires many seconds per \textsc{Ipopt} iteration to solve, while also requiring a small number of iterations (i.e., between ten and twenty).
Because of these performance characteristics, this study considers the problem to be appropriate for the proposed calibration task.

For each linear solver, the number of parallel threads was varied between two and seventy-two, in steps of two.
The executable used for solving the previously-described problem was then initiated, and the wall clock time required to reach optimality was measured.
To reduce experimental noise, each experiment was repeated five times, and the average solve time of each batch was computed.
To obtain the ``best'' number of threads for each solver, the minimum average solve time and corresponding thread count was determined.
For subsequent experiments in this study, the thread number obtained via this process is used.

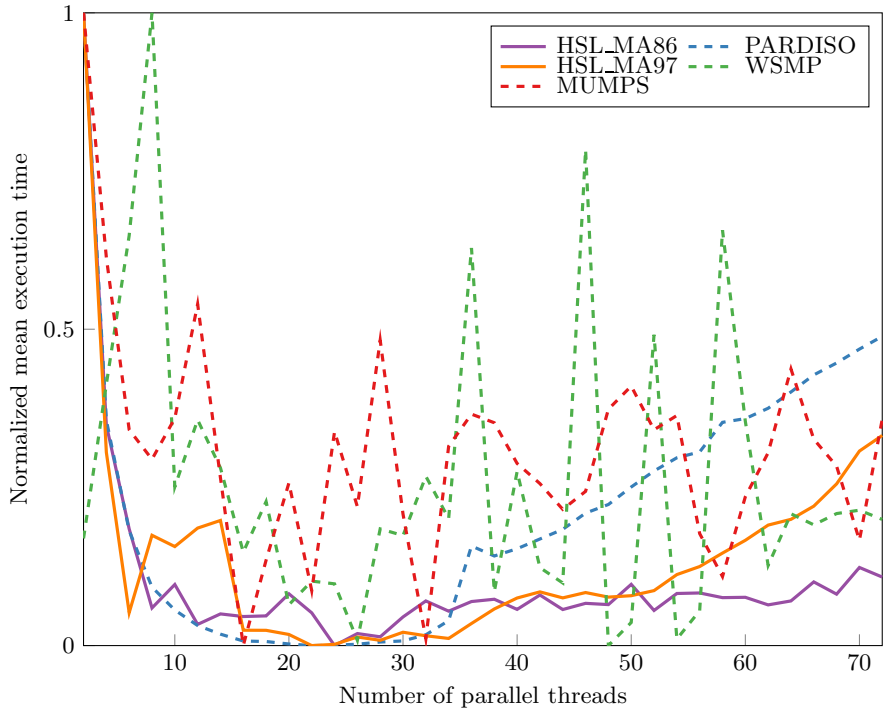
\begin{figure}[t]
    \centering
    \input{computational_experiments/figures/calibration.tex}
	 \caption{Comparison of \textsc{Ipopt} performance over various linear solvers using various numbers of parallel threads, used for the calibration in Section \ref{subsec:calibration}.}
    \label{fig:calibration-performance}
	\vspace{-1em}
\end{figure}

\begin{table}[t]
	\begin{center}
		\begin{tabular}{|c|c|c|c|c|c|}
			\hline
			\textbf{Linear Solver} & \textsc{HSL\_MA86} &
			\textsc{HSL\_MA97} & \textsc{MUMPS} & \textsc{PARDISO} & \textsc{WSMP} \\ \hline
			\textbf{Best \# Threads} & 24 & 22 & 16 & 22 & 48 \\ \hline
			\textbf{Min Solve Time (s)} & 352.1 & 749.2 & 2164.7 & 65.8 & 110.5 \\ \hline
			\textbf{Max Solve Time (s)} & 916.5 & 1197.8 & 2390.5 & 232.7 & 3423.5 \\ \hline
		\end{tabular}
		\caption{Calibrated thread counts for select parallel solvers. Here, ``Min/Max Solve Time'' correspond to experimental solve times over all tested thread counts.}
		\label{tab:linear-solver-threads}
	\end{center}
	\vspace{-3em}
\end{table}

Figure \ref{fig:calibration-performance} plots the normalized mean solve time for various numbers of parallel threads and the linear solvers considered in this section, where zero corresponds to the minimum mean solve time and one corresponds to the maximum.
Here, \textsc{HSL\_MA86} shows large times for small numbers of threads, smaller times when using around twenty-four threads, and plateauing times for larger numbers of threads.
\textsc{HSL\_MA97} shows erratic behavior for smaller numbers of threads, then generally increasing times for larger numbers.
\textsc{MUMPS} shows overall erratic behavior, and the similar minimum and maximum solve times in Table \ref{tab:linear-solver-threads} indicate that little changes while modifying \texttt{OMP\_NUM\_THREADS}.
\textsc{PARDISO} shows smooth behavior, where times generally increase when shifted further from the best twenty-two threads.
\textsc{WSMP} shows very erratic behavior, although performance is generally poorer for smaller numbers of threads.
Table \ref{tab:linear-solver-threads} summarizes the best number of threads per solver obtained via the calibration process along with minimum and maximum experimental solve times over all tested thread counts.

\subsection{Constrained and Unconstrained Testing Environment}
\label{subsec:cutest}
\citet{bongartz1995cute} introduce the Constrained and Unconstrained Testing Environment (\textsc{CUTE}), a collection of tools for testing small- and large-scale nonlinear optimization algorithms.
The initial version contains methods for comparing optimization solvers as well as $608$ different nonlinear optimization problems of various sizes and difficulty.
Furthermore, many of these problems can have their numbers of variables and/or constraints modified, allowing for more variation.
\citet{gould2003cuter} expand upon \textsc{CUTE} to develop \textsc{CUTE} (revisited; \textsc{CUTEr}), which offers backward compatibility to \textsc{CUTE}, automated multiplatform installation, interfaces to additional optimization packages (including \textsc{Ipopt}), and an annex to the original problem collection.
Finally, \citet{Gould:2015:CCU:2746412.2746437} extend \textsc{CUTEr} to develop \textsc{CUTE} (safe threads; \textsc{CUTEst}).
In addition to thread safety, the package includes a more modular design, revised installation procedure, and hundreds of new problems, totalling $1369$ instances.
These iterations of \textsc{CUTE} are often regarded as standard benchmarks for nonlinear programming.

In this study, all $1369$ \textsc{CUTEst} problem instances available as of May 2019 were solved using \textsc{Ipopt} with the \texttt{runcutest} executable.
This was repeated for each of the nine available CPU-based linear solvers, as well as the GPU-enabled version of \textsc{SPRAL}.
For each instance, default problem parameters were used, i.e., the numbers of variables and constraints were not modified from their default Standard Input Format (SIF) encodings.
Although this does not imply that all problem instances used their smallest settings, the problem set is nonetheless skewed toward easier instances, making its applicability to large-scale nonlinear optimization arguable.
Furthermore, many of the problems have more constraints than variables, for which \textsc{Ipopt} returns an exception that the problem has ``too few degrees of freedom.''
Thus, a large number of problems cannot be solved, no matter the linear solver.
Nonetheless, due to the ubiquity of \textsc{CUTEst} and the breadth of problem types contained therein, these experiments capture the general performance of linear solvers on easy to moderately difficult problems.

\begin{figure}[t]
    \centering
    \input{computational_experiments/figures/cutest.tex}
	 \caption{Comparison of \textsc{Ipopt} performance over various linear solvers using the Constrained and Unconstrained Testing Environment (safe threads; \textsc{CUTEst}).}
    \label{fig:cutest-performance}
	 \vspace{-1em}
\end{figure}
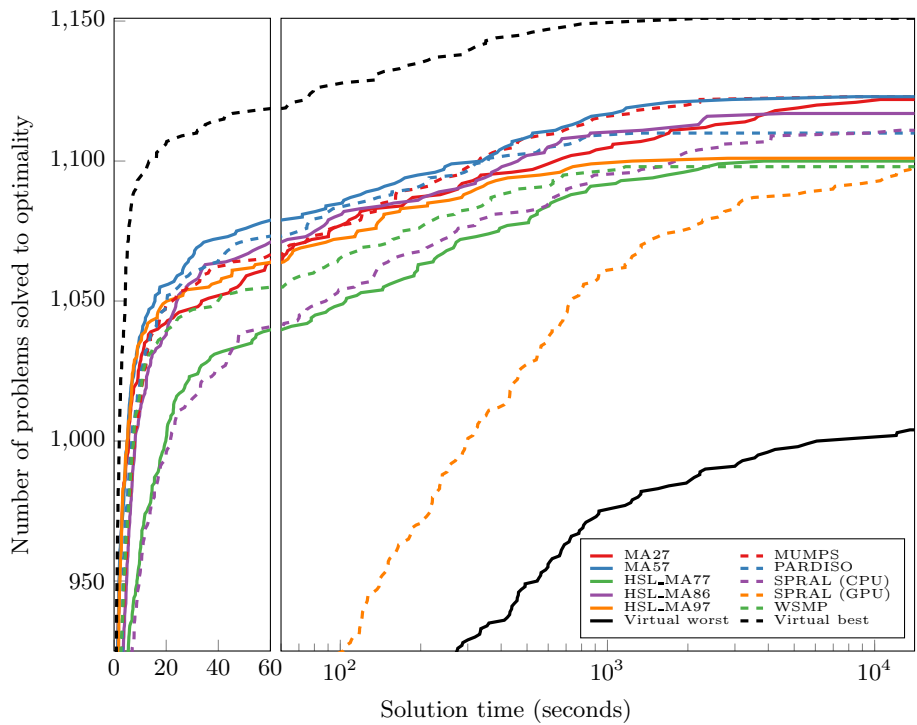

Inspired by the performance profiles presented by \citet{lundell2018supporting}, Figure \ref{fig:cutest-performance} plots the performance profiles for the \textsc{CUTEst} instances when solved with \textsc{Ipopt} and the ten compared linear solvers.
Here, the vertical axis indicates the cumulative number of instances solved to optimality on or before the time threshold presented on the horizontal axis.
Note that this method for performance profiling contrasts with the conventional method presented by \citet{dolan2002benchmarking}, which, as discussed by \citet{gould2016note} and \citet{hekmati2018nested}, can be misleading when comparing more than two profiles.
Since a large number of instances are solved prior to some small time threshold, observe that the plot is divided into two separate horizontally-scaled portions to encourage readability.
Additionally, ``virtual best'' and ``virtual worst'' profiles are included to distinguish best- and worst-case performance over all linear solvers.
This profiling method is also repeated in Sections \ref{subsec:opf}, \ref{subsec:2d-pde}, and \ref{subsec:3d-pde}.

Here, \textsc{MA57} and \textsc{MUMPS} appear to be the best performers over the majority of the problem set.
They solve large proportions of problems in small amounts of time and are also capable of solving more challenging problems.
\textsc{HSL\_MA86} and \textsc{PARDISO} perform similarly, solving large proportions of problems quickly but, overall, fewer than \textsc{MA57} and \textsc{MUMPS}.
\textsc{MA27} performs more slowly than \textsc{MA57} and \textsc{MUMPS} but eventually solves a similar number of problems.
\textsc{HSL\_MA97} and \textsc{WSMP} are generally slower and appear less robust than the aforementioned solvers.
\textsc{HSL\_MA77} performs much more slowly than the former two solvers but solves roughly the same number of problems.
\textsc{SPRAL} (CPU) also has difficulty solving easy problems quickly but eventually solves as many problems as \textsc{PARDISO}.
Finally, \textsc{SPRAL} (GPU) is the worst performer, taking long to solve even easy problems.
Presumably, the relatively good performance of the serial solvers \textsc{MA27} and \textsc{MA57} is due to the large number of easy problems, where some parallel solvers may suffer from larger overhead requirements.

\subsection{Optimal Power Flow Benchmarks}
\label{subsec:opf}
The optimization of networks constrained by nonlinear relationships or with nonlinear objectives has attracted much attention in recent decades.
Common applications are frequently seen in the optimization of potable water \citep{raghunathan2013global}, electric power \citep{coffrin2015qc}, natural gas \citep{borraz2016convex}, and traffic networks \citep{florian1976application}, to name a few.
Important applications also exist in more foundational fields such as machine learning and applied statistics, for example, in the learning of the structure and parameters associated with Ising models \citep{lokhov2018optimal}.
Network optimization problems differ from those in \textsc{CUTEst}, which encompasses a wider variety of problem classes, and PDE-constrained problems, which typically have greater topological structure.

The most common and well-studied nonlinear network optimization problems of this era are those constrained by the physics of power flow.
The most fundamental of these is the AC-OPF, which seeks to determine an optimal generation and distribution plan for power under physical constraints.
In its most comprehensive form, the AC-OPF is constrained by Ohm's and Kirchoff's laws and subject to power generation and voltage magnitude bounds at buses (nodes) and current bounds along lines (arcs).
Although the AC-OPF is nonconvex in its full form, a number of convex relaxations have been proposed to encourage tractability while maintaining solution quality.
These include the convex quadratic (QC) and second-order cone (SOC) relaxations, both of which relax the products of voltages that appear in the nonconvex AC-OPF \citep{coffrin2015qc, jabr2006radial}.

This section focuses on benchmarking convex and nonconvex OPF problems over a variety of standard network instances.
Specifically, the QC-OPF and SOC-OPF relaxations, described above, and the full nonconvex AC-OPF are benchmarked across all linear solvers and networks.
The problem formulations are implemented by \textsc{PowerModels.jl}, an open-source Julia/\textsc{JuMP} package for steady-state power network optimization \citep{coffrin2018powermodels}.
The network instances comprise version 19.05 of the ``Power Grid Lib - Optimal Power Flow'' of \citet{pglibopf}, which includes networks defined with typical, heavily-loaded, and small phase angle difference conditions.
The networks range in size from $3$ to $13,659$ buses and encompass a variety of topologies.
All experiments ($396$ per solver) were performed using Julia, \textsc{JuMP.jl}, \textsc{PowerModels.jl}, and \textsc{Ipopt.jl}.

\begin{figure}[t]
    \centering
    \input{computational_experiments/figures/opf.tex}
	 \caption{Comparison of \textsc{Ipopt} performance over various linear solvers using the ``Power Grid Lib - Optimal Power Flow'' and three different OPF formulations.}
    \label{fig:opf-performance}
	 \vspace{-2em}
\end{figure}
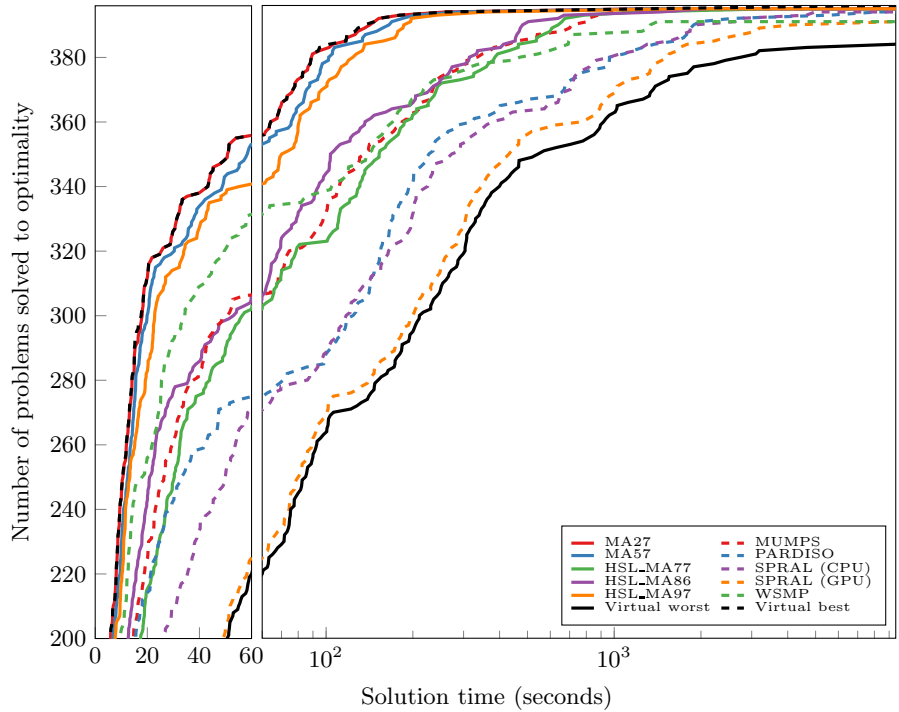

Figure \ref{fig:opf-performance} plots the performance profiles for OPF instances when solved with \textsc{Ipopt} and the ten compared linear solvers.
The behaviors of solvers are different than those observed in Figure \ref{fig:cutest-performance}.
First, \textsc{MA27} and \textsc{MA57} closely track the virtual best profile and easily solve most problems within hundreds of seconds.
\textsc{HSL\_MA97} is only somewhat slower than these two solvers.
\textsc{HSL\_MA77}, \textsc{HSL\_MA86}, \textsc{MUMPS}, and \textsc{WSMP} closely track one another and eventually solve most problem instances.
\textsc{PARDISO} and both versions of \textsc{SPRAL} are typically much slower than the aforementioned solvers.
These profiles suggest that even the largest of OPF instances result in relatively easy nonlinear programs.

\subsection{Two-dimensional Partial Differential Equations}
\label{subsec:2d-pde}
With the advent of high-performance computing, the simulation of large-scale PDEs has become ubiquitous.
In accord, there has been an increasing interest in the optimal design, control, and parameter estimation of systems governed by PDEs.
Examples of PDE-constrained optimization appear in a variety of fields, including hydrodynamics, robotics, and chemical engineering.
There are two important characteristics associated with PDE-constrained optimization problems.
The first is their large number of variables, which are required to model system state and control decisions.
The second is the sparsity of constraints, which is a result of spatial discretization via finite differencing, where each variable depends only a small number of neighboring variables.
These properties result in nonlinear programs that are typically large, highly structured, and sparse \citep{10.1007/978-3-642-55508-4_1, wachter2009short}.

This section focuses on benchmarking a number of scalable two-dimensional (2D) PDE-constrained problems that are packaged as examples within \textsc{Ipopt} under the \texttt{ScalableProblems} subdirectory.
These include eight boundary control problems adapted from \citet{maurer2000optimization} and ten distributed control problems adapted from \citet{maurer2001optimization}.
Both problem classes optimize over a set of discretized 2D elliptic PDEs.
Additionally, four examples from \citet{mittelmann2001sufficient} are included that involve 2D parabolic PDEs.
The details of the discretized NLP formulations are omitted here for brevity but are described in the aforementioned references.
For each problem, the problem dimension $N$ was scaled from $128$ to $1024$ in steps of $128$, resulting in $176$ experiments per linear solver.
For example, the largest problem, \texttt{MDistCntrl6a} with $N = 1024$, results in a nonlinear program containing $2,101,248$ variables (all of which have lower and upper bounds), $1,052,672$ equality constraints, $6,299,648$ nonzeros in the Jacobian matrix, and $3,145,728$ nonzeros in the Lagrangian Hessian matrix.

\begin{figure}[t]
    \centering
    \input{computational_experiments/figures/2d-pde.tex}
	 \caption{Comparison of \textsc{Ipopt} performance over various linear solvers using the two-dimensional partial differential equation test problem set.}
    \label{fig:2d-pde-performance}
\end{figure}
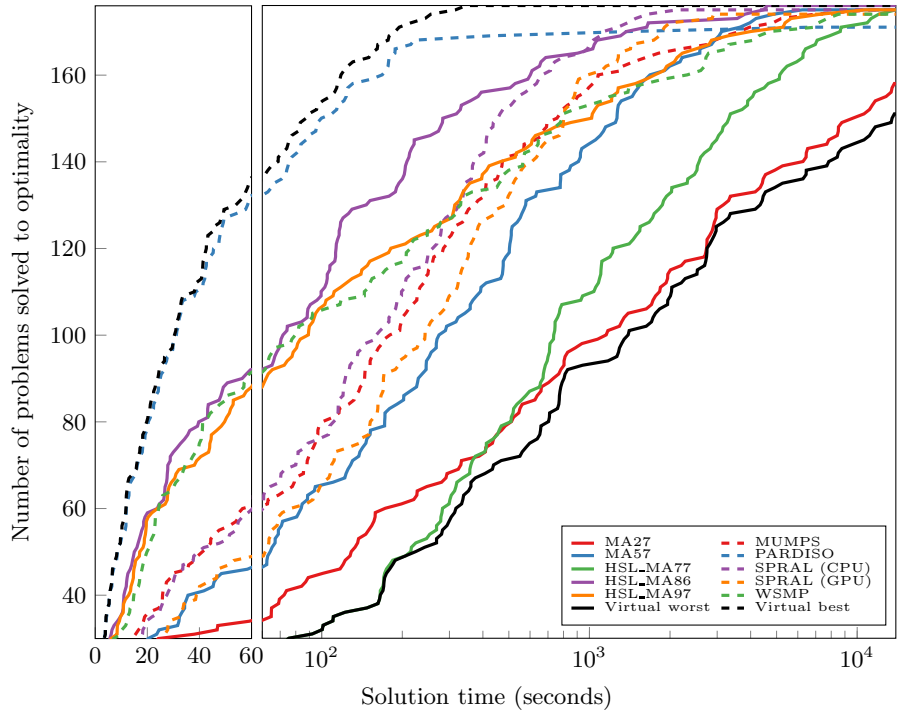

Figure \ref{fig:2d-pde-performance} plots performance profiles for the 2D PDE instances when solved with \textsc{Ipopt} and the ten linear solvers.
For most problems, \textsc{PARDISO} vastly outperforms all other solvers, but for a small number of problems, fails to converge.
\textsc{HSL\_MA86} performs well on most problems and appears more robust than \textsc{PARDISO} on some problems, solving all $176$ instances.
\textsc{HSL\_MA97} and \textsc{WSMP} solve easier problems quickly and are comparable to other solvers for most problems, solving nearly all instances.
For easier problems, \textsc{MUMPS}, \textsc{SPRAL} (CPU), and \textsc{SPRAL} (GPU) behave comparably, although \textsc{SPRAL} appears to be faster on more difficult problems.
\textsc{MA57} performs more slowly than these solvers, and \textsc{HSL\_MA77} is drastically slower. 
In this set of experiments, the drawbacks of \textsc{MA27} are finally exposed, as it generally cannot contend with other solvers on large problems.
Furthermore, the benefits of \textsc{MA57} over \textsc{MA27} are clearly seen as problem size is increased, perhaps due to the use of BLAS.

\subsection{Three-dimensional Partial Differential Equations}
\label{subsec:3d-pde}
Similar to Section \ref{subsec:2d-pde}, this section focuses on benchmarking scalable three-dimensional (3D) PDE-constrained problems that are packaged as examples within \textsc{Ipopt}.
These include four boundary control problems adapted from \citet{maurer2000optimization}, extended from 2D to 3D.
For each problem, the problem dimension $N$ was scaled from $4$ to $64$ in steps of $4$, resulting in $64$ experiments per linear solver.
The problems at their largest size (with $N = 64$) have around $290,000$ variables, $262,144$ equality constraints, $1,835,008$ or $7,075,600$ nonzeros in the Jacobian matrix, and $286,720$ or $778,240$ nonzeros in the Lagrangian Hessian matrix.
Thus, compared to the 2D PDE problems of Section \ref{subsec:2d-pde}, the problems in this section are smaller but more dense.
This generally results in greater fill-in in matrix factors, increasing the difficulty of the factorization phase.

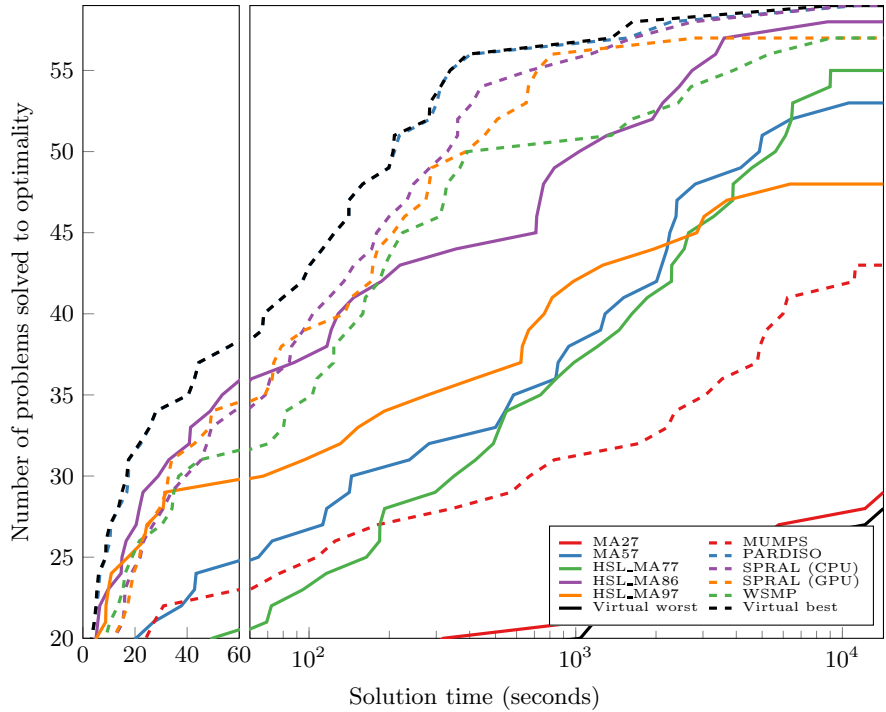
\begin{figure}[t]
    \centering
    \input{computational_experiments/figures/3d-pde.tex}
	 \caption{Comparison of \textsc{Ipopt} performance over various linear solvers using the three-dimensional partial differential equation test problem set.}
    \label{fig:3d-pde-performance}
\end{figure}

Figure \ref{fig:3d-pde-performance} plots the performance profiles for 3D PDE instances when solved with \textsc{Ipopt} and the ten linear solvers.
For most problems, \textsc{PARDISO} and the two versions of \textsc{SPRAL} are faster and more robust than other solvers.
Here, the strong performance characteristics of \textsc{PARDISO} on large problems are proven when considered with its relatively fast performance in Section \ref{subsec:2d-pde}.
For easier problems, \textsc{WSMP} and \textsc{HSL\_MA86} perform similarly and eventually solve nearly as many instances as the previous three solvers.
\textsc{HSL\_MA77} performs poorly on easy problems but shows its efficacy on more difficult problems.
Surprisingly, the serial \textsc{MA57} and parallel \textsc{HSL\_MA97} perform somewhat similarly, and \textsc{HSL\_MA97} appears less robust on more challenging problems.
\textsc{MUMPS} performs very poorly on most problems.
Finally, \textsc{MA27} performs most poorly of all and is typically capable of solving only the smallest and easiest of problems.

%% file: computational_experiments/figures/calibration.tex
\begin{tikzpicture}[]
	\centering
	\begin{axis}[legend cell align=left,enlargelimits=false,xtick pos=left,
		          ytick pos=left,height=10.0cm,ylabel=Normalized mean execution time,
	             width=\textwidth,legend style={at={(0.984, 0.98)},anchor=north east},
	             legend columns=2,xlabel=Number of parallel threads,
					 legend style={font=\small,column sep=1.5pt, row sep=-4.0pt},ymax=1.00,ytick={0.0,0.5,1.0}]
		\pgfplotstableread[col sep = comma]{computational_experiments/data/calibration/ma86.csv}{\maeightysix};
		\addplot[very thick, palette4] table [x = num_threads, y = normalized_solve_time]{\maeightysix};
		\addlegendentry{HSL\_MA86}
		\pgfplotstableread[col sep = comma]{computational_experiments/data/calibration/pardiso.csv}{\pardiso};
		\addplot[very thick, dashed, palette2] table [x = num_threads, y = normalized_solve_time]{\pardiso};
		\addlegendentry{PARDISO}
		\pgfplotstableread[col sep = comma]{computational_experiments/data/calibration/ma97.csv}{\maninetyseven};
		\addplot[very thick, palette5] table [x = num_threads, y = normalized_solve_time]{\maninetyseven};
		\addlegendentry{HSL\_MA97}
		\pgfplotstableread[col sep = comma]{computational_experiments/data/calibration/wsmp.csv}{\wsmp};
		\addplot[very thick, dashed, palette3] table [x = num_threads, y = normalized_solve_time]{\wsmp};
		\addlegendentry{WSMP}
		\pgfplotstableread[col sep = comma]{computational_experiments/data/calibration/mumps.csv}{\mumps};
		\addplot[very thick, dashed, palette1] table [x = num_threads, y = normalized_solve_time]{\mumps};
		\addlegendentry{MUMPS}
	\end{axis}
\end{tikzpicture}

%% file: computational_experiments/figures/cutest.tex
\begin{tikzpicture}[baseline]
    \begin{axis}[legend cell align=left,enlargelimits=false,xtick pos=left,
	              ytick pos=left,height=10.0cm,ylabel=Number of problems solved to optimality,
                 width=0.3\textwidth,ymin=925,ymax=1151,restrict x to domain=0:100,xmin=0.0,xmax=60.0]
        \pgfplotstableread[col sep = comma]{computational_experiments/data/cutest/ma27.csv}{\matwentyseven};
        \addplot[very thick, palette1] table [x = solve_time, y = num_problems_solved]{\matwentyseven};
        \pgfplotstableread[col sep = comma]{computational_experiments/data/cutest/ma57.csv}{\mafiftyseven};
        \addplot[very thick, palette2] table [x = solve_time, y = num_problems_solved]{\mafiftyseven};
        \pgfplotstableread[col sep = comma]{computational_experiments/data/cutest/ma77.csv}{\maseventyseven};
        \addplot[very thick, palette3] table [x = solve_time, y = num_problems_solved]{\maseventyseven};
        \pgfplotstableread[col sep = comma]{computational_experiments/data/cutest/ma86.csv}{\maeightysix};
        \addplot[very thick, palette4] table [x = solve_time, y = num_problems_solved]{\maeightysix};
        \pgfplotstableread[col sep = comma]{computational_experiments/data/cutest/ma97.csv}{\maninetyseven};
        \addplot[very thick, palette5] table [x = solve_time, y = num_problems_solved]{\maninetyseven};
        \pgfplotstableread[col sep = comma]{computational_experiments/data/cutest/mumps.csv}{\mumps};
        \addplot[very thick, dashed, palette1] table [x = solve_time, y = num_problems_solved]{\mumps};
        \pgfplotstableread[col sep = comma]{computational_experiments/data/cutest/pardiso.csv}{\pardiso};
        \addplot[very thick, dashed, palette2] table [x = solve_time, y = num_problems_solved]{\pardiso};
        \pgfplotstableread[col sep = comma]{computational_experiments/data/cutest/wsmp.csv}{\wsmp};
        \addplot[very thick, dashed, palette3] table [x = solve_time, y = num_problems_solved]{\wsmp};
        \pgfplotstableread[col sep = comma]{computational_experiments/data/cutest/spral.csv}{\spral};
        \addplot[very thick, dashed, palette4] table [x = solve_time, y = num_problems_solved]{\spral};
        \pgfplotstableread[col sep = comma]{computational_experiments/data/cutest/spralgpu.csv}{\spralgpu};
        \addplot[very thick, dashed, palette5] table [x = solve_time, y = num_problems_solved]{\spralgpu};
        \pgfplotstableread[col sep = comma]{computational_experiments/data/cutest/virtual_worst.csv}{\virtualworst};
        \addplot[very thick, solid, black] table [x = solve_time, y = num_problems_solved]{\virtualworst};
        \pgfplotstableread[col sep = comma]{computational_experiments/data/cutest/virtual_best.csv}{\virtualbest};
        \addplot[very thick, dashed, black] table [x = solve_time, y = num_problems_solved]{\virtualbest};
    \end{axis}
\end{tikzpicture}\hspace{-0.75em}
\begin{tikzpicture}[baseline]
	\begin{semilogxaxis}[legend cell align=left,enlargelimits=false,
	                     xtick pos=left,ytick pos=left,height=10.0cm,
								xlabel={\hspace{-0.2\textwidth}Solution time (seconds)},
								ytick=\empty,width=0.82\textwidth,
								xlabel absolute,legend style={at={(0.985, 0.02)},anchor=south east},
								legend image post style={scale=0.5},
								legend style={font=\tiny,column sep=1.5pt, row sep=-4.0pt},
								legend columns=2,ymin=925,ymax=1151,xmin=60.0]
        \pgfplotstableread[col sep = comma]{computational_experiments/data/cutest/ma27.csv}{\matwentyseven};
        \addplot[very thick, palette1] table [x = solve_time, y = num_problems_solved]{\matwentyseven};
        \addlegendentry{MA27}
        \pgfplotstableread[col sep = comma]{computational_experiments/data/cutest/mumps.csv}{\mumps};
        \addplot[very thick, dashed, palette1] table [x = solve_time, y = num_problems_solved]{\mumps};
        \addlegendentry{MUMPS}
        \pgfplotstableread[col sep = comma]{computational_experiments/data/cutest/ma57.csv}{\mafiftyseven};
        \addplot[very thick, palette2] table [x = solve_time, y = num_problems_solved]{\mafiftyseven};
        \addlegendentry{MA57}
        \pgfplotstableread[col sep = comma]{computational_experiments/data/cutest/pardiso.csv}{\pardiso};
        \addplot[very thick, dashed, palette2] table [x = solve_time, y = num_problems_solved]{\pardiso};
        \addlegendentry{PARDISO}
        \pgfplotstableread[col sep = comma]{computational_experiments/data/cutest/ma77.csv}{\maseventyseven};
        \addplot[very thick, palette3] table [x = solve_time, y = num_problems_solved]{\maseventyseven};
        \addlegendentry{HSL\_MA77}
        \pgfplotstableread[col sep = comma]{computational_experiments/data/cutest/spral.csv}{\spral};
        \addplot[very thick, dashed, palette4] table [x = solve_time, y = num_problems_solved]{\spral};
        \addlegendentry{SPRAL (CPU)}
        \pgfplotstableread[col sep = comma]{computational_experiments/data/cutest/ma86.csv}{\maeightysix};
        \addplot[very thick, palette4] table [x = solve_time, y = num_problems_solved]{\maeightysix};
		  \addlegendentry{HSL\_MA86}
        \pgfplotstableread[col sep = comma]{computational_experiments/data/cutest/spralgpu.csv}{\spralgpu};
        \addplot[very thick, dashed, palette5] table [x = solve_time, y = num_problems_solved]{\spralgpu};
        \addlegendentry{SPRAL (GPU)}
        \pgfplotstableread[col sep = comma]{computational_experiments/data/cutest/ma97.csv}{\maninetyseven};
        \addplot[very thick, palette5] table [x = solve_time, y = num_problems_solved]{\maninetyseven};
        \addlegendentry{HSL\_MA97}
        \pgfplotstableread[col sep = comma]{computational_experiments/data/cutest/wsmp.csv}{\wsmp};
        \addplot[very thick, dashed, palette3] table [x = solve_time, y = num_problems_solved]{\wsmp};
        \addlegendentry{WSMP}
        \pgfplotstableread[col sep = comma]{computational_experiments/data/cutest/virtual_worst.csv}{\virtualworst};
        \addplot[very thick, solid, black] table [x = solve_time, y = num_problems_solved]{\virtualworst};
        \addlegendentry{Virtual worst}
        \pgfplotstableread[col sep = comma]{computational_experiments/data/cutest/virtual_best.csv}{\virtualbest};
        \addplot[very thick, dashed, black] table [x = solve_time, y = num_problems_solved]{\virtualbest};
        \addlegendentry{Virtual best}
    \end{semilogxaxis}
\end{tikzpicture}\hfill

%% file: computational_experiments/figures/opf.tex
\begin{tikzpicture}[baseline]
    \begin{axis}[legend cell align=left,enlargelimits=false,xtick pos=left,
	              ytick pos=left,height=10.0cm,ylabel=Number of problems solved to optimality,
                 width=0.3\textwidth,ymin=200,ymax=396,restrict x to domain=0:100,xmin=0.0,xmax=60.0]
        \pgfplotstableread[col sep = comma]{computational_experiments/data/opf/ma27.csv}{\matwentyseven};
        \addplot[very thick, palette1] table [x = solve_time, y = num_problems_solved]{\matwentyseven};
        \pgfplotstableread[col sep = comma]{computational_experiments/data/opf/ma57.csv}{\mafiftyseven};
        \addplot[very thick, palette2] table [x = solve_time, y = num_problems_solved]{\mafiftyseven};
        \pgfplotstableread[col sep = comma]{computational_experiments/data/opf/ma77.csv}{\maseventyseven};
        \addplot[very thick, palette3] table [x = solve_time, y = num_problems_solved]{\maseventyseven};
        \pgfplotstableread[col sep = comma]{computational_experiments/data/opf/ma86.csv}{\maeightysix};
        \addplot[very thick, palette4] table [x = solve_time, y = num_problems_solved]{\maeightysix};
        \pgfplotstableread[col sep = comma]{computational_experiments/data/opf/ma97.csv}{\maninetyseven};
        \addplot[very thick, palette5] table [x = solve_time, y = num_problems_solved]{\maninetyseven};
        \pgfplotstableread[col sep = comma]{computational_experiments/data/opf/mumps.csv}{\mumps};
        \addplot[very thick, dashed, palette1] table [x = solve_time, y = num_problems_solved]{\mumps};
        \pgfplotstableread[col sep = comma]{computational_experiments/data/opf/pardiso.csv}{\pardiso};
        \addplot[very thick, dashed, palette2] table [x = solve_time, y = num_problems_solved]{\pardiso};
        \pgfplotstableread[col sep = comma]{computational_experiments/data/opf/wsmp.csv}{\wsmp};
        \addplot[very thick, dashed, palette3] table [x = solve_time, y = num_problems_solved]{\wsmp};
        \pgfplotstableread[col sep = comma]{computational_experiments/data/opf/spral.csv}{\spral};
        \addplot[very thick, dashed, palette4] table [x = solve_time, y = num_problems_solved]{\spral};
        \pgfplotstableread[col sep = comma]{computational_experiments/data/opf/spralgpu.csv}{\spralgpu};
        \addplot[very thick, dashed, palette5] table [x = solve_time, y = num_problems_solved]{\spralgpu};
        \pgfplotstableread[col sep = comma]{computational_experiments/data/opf/virtual_worst.csv}{\virtualworst};
        \addplot[very thick, solid, black] table [x = solve_time, y = num_problems_solved]{\virtualworst};
        \pgfplotstableread[col sep = comma]{computational_experiments/data/opf/virtual_best.csv}{\virtualbest};
        \addplot[very thick, dashed, black] table [x = solve_time, y = num_problems_solved]{\virtualbest};
    \end{axis}
\end{tikzpicture}\hspace{-0.75em}
\begin{tikzpicture}[baseline]
	\begin{semilogxaxis}[legend cell align=left,enlargelimits=false,
	                     xtick pos=left,ytick pos=left,height=10.0cm,
								xlabel={\hspace{-0.2\textwidth}Solution time (seconds)},
								ytick=\empty,width=0.82\textwidth,
								xlabel absolute,legend style={at={(0.985, 0.02)},anchor=south east},
								legend image post style={scale=0.5},
								legend style={font=\tiny,column sep=1.5pt, row sep=-4.0pt},
								legend columns=2,ymin=200,ymax=396,xmin=60.0]
        \pgfplotstableread[col sep = comma]{computational_experiments/data/opf/ma27.csv}{\matwentyseven};
        \addplot[very thick, palette1] table [x = solve_time, y = num_problems_solved]{\matwentyseven};
        \addlegendentry{MA27}
        \pgfplotstableread[col sep = comma]{computational_experiments/data/opf/mumps.csv}{\mumps};
        \addplot[very thick, dashed, palette1] table [x = solve_time, y = num_problems_solved]{\mumps};
        \addlegendentry{MUMPS}
        \pgfplotstableread[col sep = comma]{computational_experiments/data/opf/ma57.csv}{\mafiftyseven};
        \addplot[very thick, palette2] table [x = solve_time, y = num_problems_solved]{\mafiftyseven};
        \addlegendentry{MA57}
        \pgfplotstableread[col sep = comma]{computational_experiments/data/opf/pardiso.csv}{\pardiso};
        \addplot[very thick, dashed, palette2] table [x = solve_time, y = num_problems_solved]{\pardiso};
        \addlegendentry{PARDISO}
        \pgfplotstableread[col sep = comma]{computational_experiments/data/opf/ma77.csv}{\maseventyseven};
        \addplot[very thick, palette3] table [x = solve_time, y = num_problems_solved]{\maseventyseven};
        \addlegendentry{HSL\_MA77}
        \pgfplotstableread[col sep = comma]{computational_experiments/data/opf/spral.csv}{\spral};
        \addplot[very thick, dashed, palette4] table [x = solve_time, y = num_problems_solved]{\spral};
        \addlegendentry{SPRAL (CPU)}
        \pgfplotstableread[col sep = comma]{computational_experiments/data/opf/ma86.csv}{\maeightysix};
        \addplot[very thick, palette4] table [x = solve_time, y = num_problems_solved]{\maeightysix};
		  \addlegendentry{HSL\_MA86}
        \pgfplotstableread[col sep = comma]{computational_experiments/data/opf/spralgpu.csv}{\spralgpu};
        \addplot[very thick, dashed, palette5] table [x = solve_time, y = num_problems_solved]{\spralgpu};
        \addlegendentry{SPRAL (GPU)}
        \pgfplotstableread[col sep = comma]{computational_experiments/data/opf/ma97.csv}{\maninetyseven};
        \addplot[very thick, palette5] table [x = solve_time, y = num_problems_solved]{\maninetyseven};
        \addlegendentry{HSL\_MA97}
        \pgfplotstableread[col sep = comma]{computational_experiments/data/opf/wsmp.csv}{\wsmp};
        \addplot[very thick, dashed, palette3] table [x = solve_time, y = num_problems_solved]{\wsmp};
        \addlegendentry{WSMP}
        \pgfplotstableread[col sep = comma]{computational_experiments/data/opf/virtual_worst.csv}{\virtualworst};
        \addplot[very thick, solid, black] table [x = solve_time, y = num_problems_solved]{\virtualworst};
        \addlegendentry{Virtual worst}
        \pgfplotstableread[col sep = comma]{computational_experiments/data/opf/virtual_best.csv}{\virtualbest};
        \addplot[very thick, dashed, black] table [x = solve_time, y = num_problems_solved]{\virtualbest};
        \addlegendentry{Virtual best}
    \end{semilogxaxis}
\end{tikzpicture}\hfill

%% file: computational_experiments/figures/2d-pde.tex
\begin{tikzpicture}[baseline]
    \begin{axis}[legend cell align=left,enlargelimits=false,xtick pos=left,
	              ytick pos=left,height=10.0cm,ylabel=Number of problems solved to optimality,
                 width=0.3\textwidth,ymin=30,ymax=176,restrict x to domain=0:100,xmin=0.0,xmax=60.0]
        \pgfplotstableread[col sep = comma]{computational_experiments/data/2d-pde/ma27.csv}{\matwentyseven};
        \addplot[very thick, palette1] table [x = solve_time, y = num_problems_solved]{\matwentyseven};
        \pgfplotstableread[col sep = comma]{computational_experiments/data/2d-pde/ma57.csv}{\mafiftyseven};
        \addplot[very thick, palette2] table [x = solve_time, y = num_problems_solved]{\mafiftyseven};
        \pgfplotstableread[col sep = comma]{computational_experiments/data/2d-pde/ma77.csv}{\maseventyseven};
        \addplot[very thick, palette3] table [x = solve_time, y = num_problems_solved]{\maseventyseven};
        \pgfplotstableread[col sep = comma]{computational_experiments/data/2d-pde/ma86.csv}{\maeightysix};
        \addplot[very thick, palette4] table [x = solve_time, y = num_problems_solved]{\maeightysix};
        \pgfplotstableread[col sep = comma]{computational_experiments/data/2d-pde/ma97.csv}{\maninetyseven};
        \addplot[very thick, palette5] table [x = solve_time, y = num_problems_solved]{\maninetyseven};
        \pgfplotstableread[col sep = comma]{computational_experiments/data/2d-pde/mumps.csv}{\mumps};
        \addplot[very thick, dashed, palette1] table [x = solve_time, y = num_problems_solved]{\mumps};
        \pgfplotstableread[col sep = comma]{computational_experiments/data/2d-pde/pardiso.csv}{\pardiso};
        \addplot[very thick, dashed, palette2] table [x = solve_time, y = num_problems_solved]{\pardiso};
        \pgfplotstableread[col sep = comma]{computational_experiments/data/2d-pde/wsmp.csv}{\wsmp};
        \addplot[very thick, dashed, palette3] table [x = solve_time, y = num_problems_solved]{\wsmp};
        \pgfplotstableread[col sep = comma]{computational_experiments/data/2d-pde/spral.csv}{\spral};
        \addplot[very thick, dashed, palette4] table [x = solve_time, y = num_problems_solved]{\spral};
        \pgfplotstableread[col sep = comma]{computational_experiments/data/2d-pde/spralgpu.csv}{\spralgpu};
        \addplot[very thick, dashed, palette5] table [x = solve_time, y = num_problems_solved]{\spralgpu};
        \pgfplotstableread[col sep = comma]{computational_experiments/data/2d-pde/virtual_worst.csv}{\virtualworst};
        \addplot[very thick, solid, black] table [x = solve_time, y = num_problems_solved]{\virtualworst};
        \pgfplotstableread[col sep = comma]{computational_experiments/data/2d-pde/virtual_best.csv}{\virtualbest};
        \addplot[very thick, dashed, black] table [x = solve_time, y = num_problems_solved]{\virtualbest};
    \end{axis}
\end{tikzpicture}\hspace{-0.75em}
\begin{tikzpicture}[baseline]
	\begin{semilogxaxis}[legend cell align=left,enlargelimits=false,
	                     xtick pos=left,ytick pos=left,height=10.0cm,
								xlabel={\hspace{-0.2\textwidth}Solution time (seconds)},
								ytick=\empty,width=0.82\textwidth,
								xlabel absolute,legend style={at={(0.985, 0.02)},anchor=south east},
								legend image post style={scale=0.5},
								legend style={font=\tiny,column sep=1.5pt, row sep=-4.0pt},
								legend columns=2,ymin=30,ymax=176,xmin=60.0]
        \pgfplotstableread[col sep = comma]{computational_experiments/data/2d-pde/ma27.csv}{\matwentyseven};
        \addplot[very thick, palette1] table [x = solve_time, y = num_problems_solved]{\matwentyseven};
        \addlegendentry{MA27}
        \pgfplotstableread[col sep = comma]{computational_experiments/data/2d-pde/mumps.csv}{\mumps};
        \addplot[very thick, dashed, palette1] table [x = solve_time, y = num_problems_solved]{\mumps};
        \addlegendentry{MUMPS}
        \pgfplotstableread[col sep = comma]{computational_experiments/data/2d-pde/ma57.csv}{\mafiftyseven};
        \addplot[very thick, palette2] table [x = solve_time, y = num_problems_solved]{\mafiftyseven};
        \addlegendentry{MA57}
        \pgfplotstableread[col sep = comma]{computational_experiments/data/2d-pde/pardiso.csv}{\pardiso};
        \addplot[very thick, dashed, palette2] table [x = solve_time, y = num_problems_solved]{\pardiso};
        \addlegendentry{PARDISO}
        \pgfplotstableread[col sep = comma]{computational_experiments/data/2d-pde/ma77.csv}{\maseventyseven};
        \addplot[very thick, palette3] table [x = solve_time, y = num_problems_solved]{\maseventyseven};
        \addlegendentry{HSL\_MA77}
        \pgfplotstableread[col sep = comma]{computational_experiments/data/2d-pde/spral.csv}{\spral};
        \addplot[very thick, dashed, palette4] table [x = solve_time, y = num_problems_solved]{\spral};
        \addlegendentry{SPRAL (CPU)}
        \pgfplotstableread[col sep = comma]{computational_experiments/data/2d-pde/ma86.csv}{\maeightysix};
        \addplot[very thick, palette4] table [x = solve_time, y = num_problems_solved]{\maeightysix};
		  \addlegendentry{HSL\_MA86}
        \pgfplotstableread[col sep = comma]{computational_experiments/data/2d-pde/spralgpu.csv}{\spralgpu};
        \addplot[very thick, dashed, palette5] table [x = solve_time, y = num_problems_solved]{\spralgpu};
        \addlegendentry{SPRAL (GPU)}
        \pgfplotstableread[col sep = comma]{computational_experiments/data/2d-pde/ma97.csv}{\maninetyseven};
        \addplot[very thick, palette5] table [x = solve_time, y = num_problems_solved]{\maninetyseven};
        \addlegendentry{HSL\_MA97}
        \pgfplotstableread[col sep = comma]{computational_experiments/data/2d-pde/wsmp.csv}{\wsmp};
        \addplot[very thick, dashed, palette3] table [x = solve_time, y = num_problems_solved]{\wsmp};
        \addlegendentry{WSMP}
        \pgfplotstableread[col sep = comma]{computational_experiments/data/2d-pde/virtual_worst.csv}{\virtualworst};
        \addplot[very thick, solid, black] table [x = solve_time, y = num_problems_solved]{\virtualworst};
        \addlegendentry{Virtual worst}
        \pgfplotstableread[col sep = comma]{computational_experiments/data/2d-pde/virtual_best.csv}{\virtualbest};
        \addplot[very thick, dashed, black] table [x = solve_time, y = num_problems_solved]{\virtualbest};
        \addlegendentry{Virtual best}
    \end{semilogxaxis}
\end{tikzpicture}\hfill

%% file: computational_experiments/figures/3d-pde.tex
\begin{tikzpicture}[baseline]
    \begin{axis}[legend cell align=left,enlargelimits=false,xtick pos=left,
	              ytick pos=left,height=10.0cm,ylabel=Number of problems solved to optimality,
                 width=0.3\textwidth,ymin=20,ymax=59,restrict x to domain=0:100,xmin=0.0,xmax=60.0]
        \pgfplotstableread[col sep = comma]{computational_experiments/data/3d-pde/ma27.csv}{\matwentyseven};
        \addplot[very thick, palette1] table [x = solve_time, y = num_problems_solved]{\matwentyseven};
        \pgfplotstableread[col sep = comma]{computational_experiments/data/3d-pde/ma57.csv}{\mafiftyseven};
        \addplot[very thick, palette2] table [x = solve_time, y = num_problems_solved]{\mafiftyseven};
        \pgfplotstableread[col sep = comma]{computational_experiments/data/3d-pde/ma77.csv}{\maseventyseven};
        \addplot[very thick, palette3] table [x = solve_time, y = num_problems_solved]{\maseventyseven};
        \pgfplotstableread[col sep = comma]{computational_experiments/data/3d-pde/ma86.csv}{\maeightysix};
        \addplot[very thick, palette4] table [x = solve_time, y = num_problems_solved]{\maeightysix};
        \pgfplotstableread[col sep = comma]{computational_experiments/data/3d-pde/ma97.csv}{\maninetyseven};
        \addplot[very thick, palette5] table [x = solve_time, y = num_problems_solved]{\maninetyseven};
        \pgfplotstableread[col sep = comma]{computational_experiments/data/3d-pde/mumps.csv}{\mumps};
        \addplot[very thick, dashed, palette1] table [x = solve_time, y = num_problems_solved]{\mumps};
        \pgfplotstableread[col sep = comma]{computational_experiments/data/3d-pde/pardiso.csv}{\pardiso};
        \addplot[very thick, dashed, palette2] table [x = solve_time, y = num_problems_solved]{\pardiso};
        \pgfplotstableread[col sep = comma]{computational_experiments/data/3d-pde/wsmp.csv}{\wsmp};
        \addplot[very thick, dashed, palette3] table [x = solve_time, y = num_problems_solved]{\wsmp};
        \pgfplotstableread[col sep = comma]{computational_experiments/data/3d-pde/spral.csv}{\spral};
        \addplot[very thick, dashed, palette4] table [x = solve_time, y = num_problems_solved]{\spral};
        \pgfplotstableread[col sep = comma]{computational_experiments/data/3d-pde/spralgpu.csv}{\spralgpu};
        \addplot[very thick, dashed, palette5] table [x = solve_time, y = num_problems_solved]{\spralgpu};
        \pgfplotstableread[col sep = comma]{computational_experiments/data/3d-pde/virtual_worst.csv}{\virtualworst};
        \addplot[very thick, solid, black] table [x = solve_time, y = num_problems_solved]{\virtualworst};
        \pgfplotstableread[col sep = comma]{computational_experiments/data/3d-pde/virtual_best.csv}{\virtualbest};
        \addplot[very thick, dashed, black] table [x = solve_time, y = num_problems_solved]{\virtualbest};
    \end{axis}
\end{tikzpicture}\hspace{-0.75em}
\begin{tikzpicture}[baseline]
	\begin{semilogxaxis}[legend cell align=left,enlargelimits=false,
	                     xtick pos=left,ytick pos=left,height=10.0cm,
								xlabel={\hspace{-0.2\textwidth}Solution time (seconds)},
								ytick=\empty,width=0.82\textwidth,
								xlabel absolute,legend style={at={(0.985, 0.02)},anchor=south east},
								legend image post style={scale=0.5},
								legend style={font=\tiny,column sep=1.5pt, row sep=-4.0pt},
								legend columns=2,ymin=20,ymax=59,xmin=60.0]
        \pgfplotstableread[col sep = comma]{computational_experiments/data/3d-pde/ma27.csv}{\matwentyseven};
        \addplot[very thick, palette1] table [x = solve_time, y = num_problems_solved]{\matwentyseven};
        \addlegendentry{MA27}
        \pgfplotstableread[col sep = comma]{computational_experiments/data/3d-pde/mumps.csv}{\mumps};
        \addplot[very thick, dashed, palette1] table [x = solve_time, y = num_problems_solved]{\mumps};
        \addlegendentry{MUMPS}
        \pgfplotstableread[col sep = comma]{computational_experiments/data/3d-pde/ma57.csv}{\mafiftyseven};
        \addplot[very thick, palette2] table [x = solve_time, y = num_problems_solved]{\mafiftyseven};
        \addlegendentry{MA57}
        \pgfplotstableread[col sep = comma]{computational_experiments/data/3d-pde/pardiso.csv}{\pardiso};
        \addplot[very thick, dashed, palette2] table [x = solve_time, y = num_problems_solved]{\pardiso};
        \addlegendentry{PARDISO}
        \pgfplotstableread[col sep = comma]{computational_experiments/data/3d-pde/ma77.csv}{\maseventyseven};
        \addplot[very thick, palette3] table [x = solve_time, y = num_problems_solved]{\maseventyseven};
        \addlegendentry{HSL\_MA77}
        \pgfplotstableread[col sep = comma]{computational_experiments/data/3d-pde/spral.csv}{\spral};
        \addplot[very thick, dashed, palette4] table [x = solve_time, y = num_problems_solved]{\spral};
        \addlegendentry{SPRAL (CPU)}
        \pgfplotstableread[col sep = comma]{computational_experiments/data/3d-pde/ma86.csv}{\maeightysix};
        \addplot[very thick, palette4] table [x = solve_time, y = num_problems_solved]{\maeightysix};
		  \addlegendentry{HSL\_MA86}
        \pgfplotstableread[col sep = comma]{computational_experiments/data/3d-pde/spralgpu.csv}{\spralgpu};
        \addplot[very thick, dashed, palette5] table [x = solve_time, y = num_problems_solved]{\spralgpu};
        \addlegendentry{SPRAL (GPU)}
        \pgfplotstableread[col sep = comma]{computational_experiments/data/3d-pde/ma97.csv}{\maninetyseven};
        \addplot[very thick, palette5] table [x = solve_time, y = num_problems_solved]{\maninetyseven};
        \addlegendentry{HSL\_MA97}
        \pgfplotstableread[col sep = comma]{computational_experiments/data/3d-pde/wsmp.csv}{\wsmp};
        \addplot[very thick, dashed, palette3] table [x = solve_time, y = num_problems_solved]{\wsmp};
        \addlegendentry{WSMP}
        \pgfplotstableread[col sep = comma]{computational_experiments/data/3d-pde/virtual_worst.csv}{\virtualworst};
        \addplot[very thick, solid, black] table [x = solve_time, y = num_problems_solved]{\virtualworst};
        \addlegendentry{Virtual worst}
        \pgfplotstableread[col sep = comma]{computational_experiments/data/3d-pde/virtual_best.csv}{\virtualbest};
        \addplot[very thick, dashed, black] table [x = solve_time, y = num_problems_solved]{\virtualbest};
        \addlegendentry{Virtual best}
    \end{semilogxaxis}
\end{tikzpicture}\hfill

%% file: concluding_remarks/concluding_remarks.tex
\section{Concluding Remarks}
\label{sec:concluding_remarks}
This study presented a comparison of linear solvers and their performance within \textsc{Ipopt}, an interior point optimization algorithm.
Solver performance was measured over a number of problem sets, including the \textsc{CUTEst} collection; OPF problems defined by a variety of formulations and networks; and a number of scalable 2D and 3D PDE-constrained problems.
These benchmarks were intended to capture the performance of solvers over various problem classes, sizes, and difficulty.
Notably, this study is the first to compare the performance implications of selecting among the large set of linear solvers available within \textsc{Ipopt}.
In addition, the study introduced the coupling of an open-source linear solver, \textsc{SPRAL}, capable of heterogeneous parallelism over multi-core CPUs and GPUs.
Aside from \textsc{SPRAL} and \textsc{PARDISO}, options for \textsc{Ipopt} and linear solvers were left mostly at their defaults, intended to capture the behavior of typical users.

The results indicate a number of interesting characteristics, separating nonlinear programming folklore from empirical fact.
These observations motivate the solver recommendations presented in Table \ref{tab:linear-solver-suggestions}.
The first observation is that conventional nonlinear programming benchmarking suites, e.g., the \textsc{CUTEst} experiments of Section \ref{subsec:cutest} and the OPF experiments of Section \ref{subsec:opf}, do not result in matrices that are particularly large or difficult to factorize.
As such, serial solvers like \textsc{MA27} and \textsc{MA57} perform well, presumably due to the lack of overhead associated with parallel solvers.
The second observation is that, for the larger problems contained in the 2D and 3D PDE-constrained problem sets, performance sometimes differs greatly.
On the large, sparse problems contained in the 2D set, \textsc{PARDISO} performs favorably, and serial solvers (e.g., \textsc{MA57}) are surprisingly capable.
However, for more dense and difficult problems contained in the 3D set, performance varies.
Here, the multifrontal \textsc{SPRAL} also performs favorably, presumably due to its modern ordering and pivoting strategies.

Another observation is that the benefits of parallelism are sometimes suspect.
For \textsc{CUTEst} and OPF problems, many \textsc{Ipopt} iterations (tens to hundreds) are required, and the overhead of parallel solvers in reinitializing data structures could be large.
This is exaggerated for \textsc{SPRAL} (GPU), likely due to the overhead involved in shuffling data between host and GPU.
However, for difficult problems contained in the PDE-constrained sets, the benefits of parallelism are clear.
Generally, smaller numbers of \textsc{Ipopt} iterations (tens) are required, but more time is devoted to matrix factorization.
Finally, the benefits of GPU parallelism during factorization are not immediately apparent.
This could change if larger problems are considered (i.e., those requiring more than four hours to solve), but because of memory limitations, these experiments would be difficult to execute.

Based on the findings of this work, it is suggested that future work focus on creating larger and more difficult NLP benchmarks.
This may involve filtering problems from \textsc{CUTEst} and increasing the dimensions of modifiable instances.
This may also involve increasing the size of OPF benchmarks, perhaps through the use of time-extended formulations.
Future work should also focus on understanding the tradeoffs involved in GPU parallelism.
In this study, the Tesla K40M was used in benchmarks because of its prevalence on the cluster used for this work, but use of GPUs with more modern architectures may yield different results.
Finally, this study focused on the solution of problems that can be solved using shared memory.
An exploration of out-of-core techniques and distributed computation is crucial to understand performance on overly large problems.

\begin{table}[t]
	\begin{center}
		\begin{tabular}{|c|c|c|}
			\hline \textbf{Problem Type} & \textbf{Best Licensed Solver}  & \textbf{Best Free Solver} \\ \hline
			Easy & \textsc{MA57} & \textsc{MA27} \\ \hline
			Difficult & \textsc{PARDISO} & \textsc{SPRAL} \\ \hline
		\end{tabular}
		\caption{Summary of suggested linear solvers within \textsc{Ipopt} based on problem type and solver license, inferred from computational experiments in Section \ref{sec:computational_experiments}.}
		\label{tab:linear-solver-suggestions}
	\end{center}
\end{table}

%% file: acknowledgments/acknowledgments.tex
\section*{Acknowledgments}
This work was supported by the U.S. Department of Energy through the Los Alamos National Laboratory.
Los Alamos National Laboratory is operated by Triad National Security, LLC, for the National Nuclear Security Administration of the U.S. Department of Energy (Contract No. 89233218CNA000001).
Research presented in this work was supported by the Laboratory Directed Research and Development program of Los Alamos National Laboratory under project number 20170574ECR.

This work was also supported by the U.S. Department of Energy through Sandia National Laboratories.
Sandia National Laboratories is a multimission laboratory managed and operated by National Technology and Engineering Solutions of Sandia, LLC, a wholly owned subsidiary of Honeywell International, Inc., for the U.S. Department of Energy's National Nuclear Security Administration under contract DE-NA-0003525.
This paper describes objective technical results and analysis.
Any subjective views or opinions that might be expressed in the paper do not necessarily represent the views of the U.S. Department of Energy or the United States Government.

%% file: bibliography/bibliography.tex
\bibliographystyle{bibliography/splncsnat}
\bibliography{bibliography/bibliography}

%% file: appendix/appendix.tex
\appendix
\section{Appendix}
\subsection{Compilation of \textsc{SPRAL}}
The commands below describe the compilation of the \textsc{SPRAL SSIDS} linear solver.
Note that a separate set of NVCC compilation commands is provided to generate a GPU-capable library that is compatible with the linking techniques used within \textsc{Ipopt}.
\begin{lstlisting}[language=bash,showstringspaces=false,basicstyle=\tiny]
#!/bin/bash
cd ${SPRALDIR}
mkdir build && ./autogen.sh
CFLAGS=-fPIC CPPFLAGS=-fPIC CXXFLAGS=-fPIC FFLAGS=-fPIC \
FCFLAGS=-fPIC NVCCFLAGS="-shared -Xcompiler -fPIC" \
./configure --prefix="${SPRALDIR}/build" \
            --with-blas="-L${MKLROOT}/lib/intel64 -lmkl_gf_lp64 \
                         -lmkl_gnu_thread -lmkl_core -lgomp \
                         -lpthread -lm -ldl" \
            --with-metis="-L${METISDIR} -lmetis"
make && make install
nvcc -arch=sm_35 -dlink -Xcompiler -fPIC \
    src/blas_iface.o src/core_analyse.o src/compat.o src/cuda/api_wrappers.o \
    src/cuda/cuda.o src/hw_topology/hw_topology.o \
    src/hw_topology/guess_topology.o src/lapack_iface.o src/lsmr.o \
    interfaces/C/lsmr.o src/match_order.o \
    src/matrix_util.o interfaces/C/matrix_util.o src/metis4_wrapper.o \
    src/omp.o src/pgm.o src/rutherford_boeing.o \
    interfaces/C/rutherford_boeing.o src/scaling.o interfaces/C/scaling.o \
    src/timer.o src/random.o interfaces/C/random.o src/random_matrix.o \
    interfaces/C/random_matrix.o src/ssids/gpu/alloc.o \
    src/ssids/gpu/cpu_solve.o src/ssids/gpu/datatypes.o \
    src/ssids/gpu/dense_factor.o src/ssids/gpu/factor.o \
    src/ssids/gpu/interfaces.o src/ssids/gpu/solve.o src/ssids/gpu/smalloc.o \
    src/ssids/gpu/subtree.o src/ssids/gpu/kernels/assemble.o \
    src/ssids/gpu/kernels/dense_factor.o src/ssids/gpu/kernels/reorder.o \
    src/ssids/gpu/kernels/solve.o src/ssids/gpu/kernels/syrk.o \
    src/ssids/akeep.o src/ssids/anal.o src/ssids/contrib.o \
    src/ssids/contrib_free.o src/ssids/datatypes.o src/ssids/fkeep.o \
    src/ssids/inform.o src/ssids/profile.o src/ssids/profile_iface.o \
    src/ssids/ssids.o src/ssids/subtree.o src/ssids/cpu/cpu_iface.o \
    src/ssids/cpu/NumericSubtree.o src/ssids/cpu/subtree.o \
    src/ssids/cpu/SymbolicSubtree.o src/ssids/cpu/ThreadStats.o \
    src/ssids/cpu/kernels/cholesky.o src/ssids/cpu/kernels/ldlt_app.o \
    src/ssids/cpu/kernels/ldlt_nopiv.o src/ssids/cpu/kernels/ldlt_tpp.o \
    src/ssids/cpu/kernels/wrappers.o interfaces/C/ssids.o src/ssmfe/ssmfe.o \
    src/ssmfe/core.o src/ssmfe/expert.o interfaces/C/ssmfe.o \
    interfaces/C/ssmfe_core.o interfaces/C/ssmfe_expert.o -o dlink.o
nvcc -Xcompiler -fPIC -arch=sm_35 --lib dlink.o src/blas_iface.o \
    src/core_analyse.o src/compat.o src/cuda/api_wrappers.o src/cuda/cuda.o \
    src/hw_topology/hw_topology.o src/hw_topology/guess_topology.o \
    src/lapack_iface.o src/lsmr.o interfaces/C/lsmr.o src/match_order.o \
    src/matrix_util.o interfaces/C/matrix_util.o src/metis4_wrapper.o src/omp.o \
    src/pgm.o src/rutherford_boeing.o interfaces/C/rutherford_boeing.o \
    src/scaling.o interfaces/C/scaling.o src/timer.o src/random.o \
    interfaces/C/random.o src/random_matrix.o interfaces/C/random_matrix.o \
    src/ssids/gpu/alloc.o src/ssids/gpu/cpu_solve.o src/ssids/gpu/datatypes.o \
    src/ssids/gpu/dense_factor.o src/ssids/gpu/factor.o \
    src/ssids/gpu/interfaces.o src/ssids/gpu/solve.o src/ssids/gpu/smalloc.o \
    src/ssids/gpu/subtree.o src/ssids/gpu/kernels/assemble.o \
    src/ssids/gpu/kernels/dense_factor.o src/ssids/gpu/kernels/reorder.o \
    src/ssids/gpu/kernels/solve.o src/ssids/gpu/kernels/syrk.o \
    src/ssids/akeep.o src/ssids/anal.o src/ssids/contrib.o \
    src/ssids/contrib_free.o src/ssids/datatypes.o src/ssids/fkeep.o \
    src/ssids/inform.o src/ssids/profile.o src/ssids/profile_iface.o \
    src/ssids/ssids.o src/ssids/subtree.o src/ssids/cpu/cpu_iface.o \
    src/ssids/cpu/NumericSubtree.o src/ssids/cpu/subtree.o \
    src/ssids/cpu/SymbolicSubtree.o src/ssids/cpu/ThreadStats.o \
    src/ssids/cpu/kernels/cholesky.o src/ssids/cpu/kernels/ldlt_app.o \
    src/ssids/cpu/kernels/ldlt_nopiv.o src/ssids/cpu/kernels/ldlt_tpp.o \
    src/ssids/cpu/kernels/wrappers.o interfaces/C/ssids.o src/ssmfe/ssmfe.o \
    src/ssmfe/core.o src/ssmfe/expert.o interfaces/C/ssmfe.o \
    interfaces/C/ssmfe_core.o interfaces/C/ssmfe_expert.o -o libspralgpu.a
\end{lstlisting}

\subsection{Compilation of \textsc{Ipopt} with Sequential BLAS}
The commands below describe the compilation of \textsc{Ipopt} when using sequential variants of the Intel MKL BLAS routines.
\begin{lstlisting}[language=bash,showstringspaces=false,basicstyle=\tiny]
export MKLSEQ="-L${MKLROOT}/lib/intel64 -Wl,--no-as-needed -lmkl_gf_lp64 \
               -lmkl_sequential -lmkl_core -lpthread -lm -ldl"
cd ${IPOPTSEQDIR}
ADD_CFLAGS="-fopenmp -fPIC -fno-common -fexceptions" \
ADD_CXXFLAGS="-fopenmp -fPIC -fno-common -fexceptions" \
ADD_FFLAGS="-fopenmp -fPIC -fexceptions" \
../configure --with-pic \
             --prefix="${IPOPTSEQDIR}" \
             --with-pardiso="-L${PARDISODIR} -lpardiso600-GNU720-X86-64 \
                             -lgfortran -lquadmath" \
             --with-wsmp="${WSMPDIR}/libwsmp64.a" \
             --with-blas-lib="${MKLSEQ} -lmkl_lapack95_lp64" \
             --with-lapack-lib="${MKLSEQ} -lmkl_lapack95_lp64"
make -j $(nproc --all) && make install
cd ${IPOPTSEQDIR}/Ipopt/examples/ScalableProblems && make
\end{lstlisting}

\subsection{Compilation of \textsc{Ipopt} with Parallel BLAS}
The commands below describe the compilation of \textsc{Ipopt} when using parallel variants of the Intel MKL BLAS routines.
\begin{lstlisting}[language=bash,showstringspaces=false,basicstyle=\tiny]
export MKLPAR="-L${MKLROOT}/lib/intel64 -Wl,--no-as-needed -lmkl_gf_lp64 \
               -lmkl_gnu_thread -lmkl_core -lgomp -lpthread -lm -ldl"
ADD_CFLAGS="-fopenmp -fPIC -fno-common -fexceptions" \
ADD_CXXFLAGS="-fopenmp -fPIC -fno-common -fexceptions" \
ADD_FFLAGS="-fopenmp -fPIC -fexceptions" \
../configure --with-pic \
             --prefix="${IPOPTPARDIR}" \
             --with-spral="-L${SPRALDIR} -L${CUDA_HOME}/lib64 -L${METISDIR} \
                           -L${HWLOCDIR}/build/lib -lspralgpu -lmetis \
                           -lcudadevrt -lcudart -lcuda -lcublas -lhwloc \
                           -lgfortran -lOpenCL" \
             --with-blas-lib="${MKLPAR} -lmkl_lapack95_lp64" \
             --with-lapack-lib="${MKLPAR} -lmkl_lapack95_lp64"
make -j $(nproc --all) && make install
cd ${IPOPTPARDIR}/Ipopt/examples/ScalableProblems && make
\end{lstlisting}